\documentclass[leqno]{amsart}

\usepackage{amscd,amssymb}
\usepackage{verbatim}
\renewcommand{\baselinestretch}{1.2}
\newtheorem{Lemma}{Lemma}[section]
\newtheorem{Th}[Lemma]{Theorem}
\newtheorem{Prop}[Lemma]{Proposition}

\newtheorem{Cor}[Lemma]{Corollary}
\newtheorem{Def}[Lemma]{Definition}
\newtheorem{Ex}[Lemma]{Example}
\newtheorem{Exs}[Lemma]{Examples}

\newtheorem{Remark}[Lemma]{Remark}
\newenvironment{Proof}{{\sc Proof.}\ }{~\rule{1ex}{1ex}\vspace{0.5truecm}}

\newcommand{\End}{\mbox{\rm End}}

\newcommand{\Hom}{\mbox{\rm Hom}}
\newcommand{\Ext}{\mbox{\rm Ext}}
\newcommand{\Tor}{\mbox{\rm Tor}}

\newcommand{\add}{\mbox{\rm add}}

\newcommand{\Add}{\mbox{\rm Add}}
\newcommand{\Prod}{\mbox{\rm Prod}}

\newcommand{\B}{\mathcal{B}}

\newcommand{\coker}{\mbox{\rm coker}}
\newcommand{\Ker}{\mathrm{ker}}

\newcommand{\im}{\mbox{\rm Im}}

\newcommand{\N}{\mathbb N}
\newcommand{\Z}{\mathbb{Z}}
\newcommand{\C}{\mathcal{C}}
\newcommand{\D}{\mathcal{D}}

\newcommand{\Q}{\mathcal{Q}}
\newcommand{\A}{\mathcal{A}}

\newcommand{\Mod}{\mbox{\rm Mod-}}
\newcommand{\lmod}{\mbox{\rm -mod}}

\newcommand{\LMod}{\mbox{\rm -Mod}}
\newcommand{\rmod}{\mbox{\rm mod-}}

\newcommand{\Fcal}{\ensuremath{\mathcal{F}}}
\newcommand{\Ecal}{\ensuremath{\mathcal{E}}}

\newcommand{\Scal}{\ensuremath{\mathcal{S}}}

\newcommand{\Tcal}{\ensuremath{\mathcal{T}}}
\newcommand{\Mcal}{\ensuremath{\mathcal{M}}}
\newcommand{\Lcal}{\ensuremath{\mathcal{L}}}
\newcommand{\Ccal}{\ensuremath{\mathcal{C}}}
\newcommand{\Qcal}{\ensuremath{\mathcal{Q}}}
\newcommand{\p}{\ensuremath{\mathbf{p}}}
\newcommand{\q}{\ensuremath{\mathbf{q}}}
\newcommand{\tube}{\ensuremath{\mathbf{t}}}

\newcommand{\ra}{\rightarrow}
\newcommand{\lora}{\longrightarrow}
\newcommand{\mapr}[1]{\stackrel{#1}{\longrightarrow}}

\newcommand{\exs}[5]{0\longrightarrow #1 \mapr{#2} #3 \mapr{#4} #5\longrightarrow 0}

\newcommand{\M}{\mbox{Mod\,$\Lambda$}}

\newcommand{\eme}{{\mathcal M}}

\begin{document}

\title{Mittag-Leffler conditions on modules}

\author{Lidia Angeleri H\"ugel}
\address{Dipartimento di Informatica e
Comunicazione, Universit\`a degli Studi dell'Insubria, Via Mazzini
5, I - 21100 Varese, Italy} \email{lidia.angeleri@uninsubria.it}
\author{Dolors Herbera}
\address{ Departament de Matem\`atiques \\
Universitat Aut\`onoma de Barcelona \\ E-08193 Be\-lla\-te\-rra
(Barcelona), Spain} \email{dolors@mat.uab.es}

\thanks{The research of this paper was initiated while the second author was spending a sabbatical year
at the Universit\`a di Padova (Italy). She thanks her host for its
kind hospitality. She also thanks  the Universit\`a dell'Insubria, Varese
(Italy) for its hospitality in several visits while the paper was
being written. The research was continued during the visit
of the authors to CRM Barcelona in September 2006 supported by the
Research Programme on Discrete and Continuous Methods of Ring
Theory. First author partially supported by Universit\`a di Padova,
Progetto di Ateneo CDPA048343, and by PRIN 2005 "Prospettive in
teoria degli anelli, algebre di Hopf e categorie di moduli". Both
authors partially supported by the DGI and the European Regional
Development Fund, jointly, through Project
 MTM2005--00934, and by the Comissionat per Universitats i Recerca
of the Generalitat de Ca\-ta\-lunya, Project 2005SGR00206. 
\newline 2000
Mathematics Subject Classification: 16D70, 16L30, 18E15.}

\date{\today}

\begin{abstract}
We study Mittag-Leffler conditions on modules providing relative versions of classical results by Raynaud and Gruson.
We then apply our investigations to several contexts. First of all,   we give a new argument for solving the
Baer splitting problem. Moreover, we show that modules arising in cotorsion pairs satisfy certain Mittag-Leffler
conditions. In particular, this implies that tilting modules satisfy a useful finiteness condition over
their endomorphism ring. In the final section, we focus on a special tilting cotorsion pair related to the
pure-semisimplicity conjecture.

\end{abstract}

\maketitle

\renewcommand{\baselinestretch}{0.9}
 \small

\noindent{\bf Contents:}

\begin{tabbing}
{\rm 1.}\; \={$\mathcal {Q}$-Mittag-Leffler modules}\hspace{4.5cm}\=\pageref{ML}\\
{\rm 2.}\> {$H$-subgroups}\>\pageref{H}\\
{\rm 3.}\> {$\mathcal {B}$-stationary modules}\>\pageref{stat}\\
{\rm 4.} \>{Dominating maps}\>\pageref{dom}\\
{\rm 5. }\>{$\mathcal {Q}$-Mittag-Leffler modules revisited}\>\pageref{revisit}\\
{\rm 6.}\> {Relating $\mathcal {B}$-stationary and $\mathcal {Q}$-Mittag-Leffler modules}\>\pageref{relat}\\
{\rm 7.}\> {Baer modules}\>\pageref{baer}\\
{\rm 8.}\> {Matrix subgroups}\>\pageref{smatrix}\\
{\rm 9.}\> {Cotorsion pairs}\>\pageref{cotorsion}\\
{\rm 10.}\> {The pure-semisimplicity conjecture}\>\pageref{pss}
\end{tabbing}


\normalsize

\renewcommand{\baselinestretch}{1.2}

\normalsize

\section*{Introduction}\label{intro}

In the last few years, Mittag-Leffler conditions on modules were successfully employed in a number of different problems
ranging from tilting theory
to commutative algebra, and to a conjecture originating in algebraic topology. Indeed, the translation of certain
homological properties of modules into Mittag-Leffler conditions was a key step in solving the Baer splitting
problem raised by Kaplansky in 1962 \cite{abh}, as well as in proving that every  tilting class is determined
by a class of finitely presented modules \cite{bazher} \cite{BaSt}, and it is part of the strategy for
tackling the telescope conjecture for module categories  \cite{tel}.

Motivated by these results, in this paper we undertake a systematic study
of such conditions, and we give further applications of these tools.
In fact, we give a new proof for the result in \cite{abh}. Moreover, using the theory of matrix subgroups,
we provide a new interpretation of certain finiteness conditions of a module over its endomorphism ring,
in particular of endofiniteness. Furthermore, we show that Mittag-Leffler conditions
appear naturally  in the theory of cotorsion
pairs, that is, pairs of classes of modules that are orthogonal to each other with respect to the Ext functor.
As a consequence, we discover a new finiteness condition satisfied by tilting modules. Finally, we employ our
investigations to discuss the
pure-semisimplicity conjecture, developing an idea from   \cite{key}.

Further applications of   our work to finite-dimensional hereditary algebras, and to  cotorsion pairs given by
modules of bounded projective
dimension   will appear in \cite{baerml} and \cite{bazher2}, respectively.

\medskip

We now give some details on the conditions we are going to investigate.
Raynaud and Gruson studied in \cite{RG} the right modules $M$ over
a ring $R$ having the property that the canonical map $$\rho
\colon M\bigotimes_R \prod _{i\in I}Q_i\to \prod _{i\in
I}(M\bigotimes_RQ_i)$$ is injective for any family of left
$R$-modules $\{Q_i\}_{i\in I}$. They showed that this is the case
if and only if $M$ is the direct limit of a direct system
$(F_\alpha ,f_{\beta \, \alpha})_{\beta \, \alpha \in \Lambda}$ of
finitely presented modules such that
 the inverse system \[(\Hom_R(F_\alpha, B),\Hom_R(f_{\beta \, \alpha}, B))_{\beta \, \alpha \in \Lambda}\]
satisfies the Mittag-Leffler condition   for
any right $R$-module $B$. Therefore   such   modules $M$ are said to be Mittag-Leffler modules.

In this paper, we study relative versions of these properties by restricting the choice of the family
  $\{Q_i\}_{i\in I}$ and of $B$. We thus consider the notions of a
\emph{${\Q}$-Mittag-Leffler module} and of a \emph{$\B$-stationary module}. Part of our work consists in
developing these notions following closely
  \cite{RG}.

While the definition of a $\Q$-Mittag-Leffler module relies on the injectivity of the natural transformation $\rho$,
the $\mathcal{B}$-stationary modules are not ``canonically" defined. We introduce the stronger
notion of \emph{strict $\mathcal{B}$-stationary modules}. Again, this notion is inspired by \cite{RG}. Indeed, if $\mathcal{B}$ is the class of all right modules, then the strict  $\mathcal{B}$-stationary modules
are precisely the strict Mittag-Leffler modules introduced by Raynaud and Gruson, and later  studied by Azumaya \cite{Az3} and  other authors  under the name of
locally projective modules. We characterize  strict $\mathcal{B}$-stationary modules  in terms of
the injectivity of the natural transformation
\[\nu =\nu {(M,B,V)} \colon M\otimes _R\mathrm{Hom}_S(B,V)\to \mathrm{Hom}_S(\mathrm{Hom}_R(M,B),V).\]
This relates our investigations to results on matrix subgroups obtained by Zimmermann in \cite{Z}.

As mentioned above,
our original motivation are
the results in \cite{bazher}, where it was made apparent that for  a countably presented module
$M$ the vanishing of $\mathrm{Ext}_R^1(M,B)$ for all modules $B$ belonging to a class   $\mathcal{B}$ closed
under direct sums, can be characterized in terms of   $\mathcal{B}$-stationarity, see Theorem~\ref{BH} for
a precise statement. Furthermore, also the vanishing   of $\varprojlim ^1$, the first derived functor of the
inverse limit $\varprojlim$, can be interpreted in this way, see \cite{e} and \cite{abh}. We believe that a thorough understanding of
$\mathcal{B}$-stationary modules and of their relationship with strict $\mathcal{B}$-stationary  and $\Q$-Mittag-Leffler modules will provide a new insight in
 problems related to the vanishing of certain homological functors. The applications we present in this paper are oriented towards such developments.

 Let us illustrate such applications by focussing on cotorsion pairs. Let $\Scal$ be a set of finitely presented modules, and let $(\Mcal, \Lcal)$ be the cotorsion pair generated by $\Scal$. In other words, $\Lcal$ is the class of modules defined by the vanishing of $\mathrm{Ext}_R^1(\Scal,-)$, while $\Mcal$ is defined by the vanishing of $\mathrm{Ext}_R^1(-,\Lcal)$, see Definition \ref{cotorsionpair}. Denote further by $\C$ the class defined by the vanishing of $\mathrm{Tor}_1^R(\Scal,-)$. We prove in Theorem \ref{cotpair} that a module is $\Lcal$-stationary if and only if it is $\C$-Mittag-Leffler. Moreover, it turns out that every module in $\Mcal$ is strict $\Lcal$-stationary.

 In particular, this applies to cotorsion pairs arising in tilting theory (Corollary~\ref{tilting}), yielding that every tilting module $T$ is strict $T$-stationary. The latter property can be interpreted in terms of matrix subgroups and allows us to show that certain tilting modules are noetherian over their endomorphism ring, see Proposition \ref{know} and \cite{baerml}.

\medskip

The paper is organized as follows. In Section~\ref{ML} we introduce $\mathcal{Q}$-Mittag-Leffler modules, and we study the closure
properties of the class $\mathcal{Q}$ and of the class of $\mathcal{Q}$-Mittag-Leffler modules.
For our applications to cotorsion pairs, it is relevant to note the good behavior of
Mittag-Leffler modules with respect to filtrations
established
in Proposition~\ref{MLfilters}. We revisit the topic of $\mathcal{Q}$-Mittag-Leffler modules in
Section~\ref{revisit}, where we characterize them in the spirit of \cite{RG}: since the map
$\rho$ is bijective  when $M$ is finitely presented, and since every module  is a direct limit of
finitely presented modules,   one has to determine the ``gluing'' conditions on the canonical  maps $u_\alpha, u_{\beta \,\alpha}$ in the direct limit presentation  $(M, (u_\alpha)_{\alpha \in I})=\varinjlim(F_\alpha,u_{\beta \,\alpha})_{\beta \,\alpha \in I  }$
 that imply the injectivity of $\rho$. These conditions  are called
 \emph{dominating with respect to $\mathcal{Q}$}. We introduce them in Section~\ref{dom} where we also
study their basic properties.

  $\mathcal{B}$-stationary modules are introduced in Section~\ref{stat}. Hereby we  adopt the language of $H$-subgroups from \cite{zbcn}, which is the topic of Section~\ref{H}.  Our first
aim is to give an intrinsic characterization of $\mathcal{B}$-stationarity. This characterization, obtained  in Theorem~\ref{2.1.4}, is  also given in terms of   \emph{dominating maps}. It will allow  us to study the
interplay between the conditions $\mathcal{B}$-stationary
  and   $\mathcal{Q}$-Mittag-Leffler in
Section~\ref{relat}.

The interrelationship between the different conditions is further pursued in Section~\ref{cotorsion}, after having
introduced and characterized the strict $\mathcal{B}$-stationary modules in
Section~\ref{smatrix}. Note that the condition strict $\mathcal{B}$-stationary has again a good behavior under filtrations, cf.
Proposition~\ref{strictfilters}.   This intertwines our investigations
 with the theory  of cotorsion pairs. Our main
results in this context are Theorem~\ref{cotpair} and its application to tilting cotorsion pairs   in
Corollary~\ref{tilting}, which we have already described above.

A further important application is given
in
Section~\ref{baer} which is devoted to Baer modules over domains.
 A module $M$ over a commutative domain $R$ is said to be a Baer module if $\mathrm{Ext}_R^1(M,T)=0$ for any torsion
  module $T$. Kaplansky in \cite{K} raised the question whether Baer modules are projective. The last step in  the positive solution of Kaplansky's
  problem  was made in \cite{abh}. In the present work, we prove that a
countably generated Baer module over an arbitrary commutative domain is always a  Mittag-Leffler module. This
yields another proof  of the fact that Baer modules over  commutative domains are projective.

Let us mention that the techniques introduced by Raynaud and Gruson
have also been used by Drinfeld in \cite{drinfeld}. We give in
Corollary~\ref{cp} a detailed proof of \cite[Theorem~2.2]{drinfeld}.

Finally, as a last application, we consider  left pure-semisimple hereditary rings in
Section~\ref{pss}. In particular,  we use
Corollary~\ref{tilting} to study the tilting cotorsion pair generated by the preprojective right modules, following an idea from \cite{key}.

\bigskip

{\bf Acknowledgements.} We would like to thank Javier S\'anchez
Serd\`a for   valuable discussions on the paper \cite{RG}, and
Silvana Bazzoni for many comments on  preliminary versions of the
paper.

\bigskip

{\bf Notation.}
 Let $R$ be a ring. Denote by $\Mod R$ the category of all right
$R$-modules, and by $\rmod R$ the subcategory of all modules
possessing a projective resolution consisting of finitely generated modules.
$R\LMod$ and   $R\lmod$ are defined correspondingly.

\medskip

For a right $R$-module $M$, we denote by
$M^\ast=\Hom_\Z(M,\mathbb{Q}/\Z)$ its character module. Instead of
the character module we can also consider another dual module, for
example, for modules  over an artin algebra $\Lambda$ we can take
$M^\ast=D(M)$ where  $D$ denotes the usual duality. If $\Scal$ is a
class of modules, we denote by $\Scal^\ast$ the corresponding class
of all duals $B^\ast$ of modules $B\in\Scal$.

\medskip

For a class $\Mcal\subset \Mod R$ we set
\begin{align*}
&\Mcal^\circ=\{X\in\Mod R\,|\, \Hom_R(M,X)=0 \, \text{ for all } M\in \Mcal\}\\
 &\,^{\circ}\Mcal=\{X\in\Mod R\,|\, \Hom_R(X,M)=0\, \text{ for all } M\in \Mcal\}\\
  &\Mcal^\perp=\{X\in\Mod R\,|\, \Ext_R^1(M,X)=0 \, \text{ for all } M\in \Mcal\}\\
 &\,^{\perp}\Mcal=\{X\in\Mod R\,|\, \Ext_R^1(X,M)=0\, \text{ for all } M\in \Mcal\}\\
  &\Mcal^\intercal=\{X\in\Mod R\,|\, \Tor^R_1(M,X)=0 \, \text{ for all } M\in \Mcal\}\\
 &\,^{\intercal}\Mcal=\{X\in\Mod R\,|\, \Tor^R_1(X,M)=0\, \text{ for all } M\in \Mcal\}
\end{align*}
Moreover, we denote by
{Add}\,{$\Mcal$}
(respectively,  {add}\,{$\Mcal$})
the class consisting of all modules isomorphic to direct summands of (finite)
direct sums of modules of ${\Mcal} $. The class  consisting of   all modules isomorphic to direct  summands of products of modules of ${\Mcal} $
is denoted by
{Prod}\,{$\Mcal$}.
Further, Gen{$\Mcal$} and Cogen{$\Mcal$} denote the class of modules generated, respectively cogenerated, by modules of ${\Mcal} $.
If $\Mcal$ consists of a single module $M$, we just write $M^\perp$, Add$M$, Prod$M$, etc.
Finally, we write ${\varinjlim
\mathcal M}$ for the class of all modules $D$ such that
$D = \varinjlim _{i \in I} M_i$ where
$\{ M_i \mid i \in I \}$
is a direct system of modules from $\mathcal M$.

\medskip

We will say that a module $M_R$ with endomorphism ring $S$ is {\em endonoetherian} if $M$ is noetherian when viewed as a left $S$-module.
If $_SM$ has finite length, then we say that $M$ is  {\em endofinite}.

\section{$\Q$-Mittag-Leffler modules} \label{ML}

\begin{Def} {\rm \cite{rothmaler} Let $M$ be a right module over a ring $R$, and let
$\Q$ be a class of left $R$-modules. We say that
$M$ is a \emph{$\Q$-Mittag-Leffler module} if the
canonical map
\[\rho \colon  M\bigotimes_R \prod _{i\in I}Q_i\to \prod _{i\in
I}(M\bigotimes_RQ_i)\] is injective for any family $\{Q_i\}_{i\in
I}$ of modules in $\Q$.
If $\Q$ just consists of a single module $Q$, then we say that $M$ is $Q$-Mittag-Leffler.
}
\end{Def}

We will need the following Lemma.

\begin{Lemma}\label{zeros} Let $M_R$ and ${}_RQ$ be a right and a
left $R$-module, respectively. Assume that $Q=\varinjlim
(K_\alpha, f_{\beta \, \alpha} )_{\beta \, \alpha\in I  }$. For
any $\alpha \in I  $, let $f_\alpha \colon  K_\alpha \to Q$ be the
induced map.

Let $q_1,\dots ,q_n\in Q$ and $x_1,\dots ,x_n\in M$ be such that
$\sum _{i=1}^nx_i\otimes q_i$ is the zero element of $M\otimes
_RQ$. Then there exist $\alpha_0 \in I  $ and $k_1,\dots ,k_n
\in K_{\alpha_0}$ such that $\sum _{i=1}^nx_i\otimes k_i$ is the
zero element of $M\otimes _RK_{\alpha_0}$ and
$f_{\alpha_0}(k_i)=q_i$ for every $i=1,\dots ,n$.
\end{Lemma}
\begin{Proof} Choose $\beta$ such that $\{q_1,\dots ,q_n\}\subseteq f_\beta (K_\beta
)$. For every $i\in \{1,\dots ,n\}$, let $k'_i\in K_\beta$ be such
that $f_\beta (k'_i)=q_i$. Since $\sum _{i=1}^nx_i\otimes q_i=0$
in $M\otimes _RQ\cong \varinjlim (M\otimes _RK_\alpha)$, there
exists ${\alpha_0} \ge \beta$ such that
\[0=(M\otimes f_{{\alpha_0} \, \beta})(\sum _{i=1}^nx_i\otimes k'_i)=\sum _{i=1}^nx_i\otimes f_{{\alpha_0} \, \beta}(k'_i).\]
Now ${\alpha_0}$ and $k_i=f_{{\alpha_0} \, \beta}(k'_i)$,
$i=1,\dots ,n$, satisfy the desired properties.
\end{Proof}

Here are some closure properties of the class $\Q$. Related results can be found
in work of Rothmaler \cite[Theorem~2.2, Remark~2.3]{rothmaler} and Zimmermann \cite[2.2]{Z}.

\begin{Th} \label{definable}
Let $R$ be a ring,   and $\Q\subseteq R\LMod$. Assume that $M\in\Mod R$ is $\Q$-Mittag-Leffler.
 Then the following statements hold true.
\begin{itemize}
\item[(i)]  $M$ is
$\Q'$-Mittag-Leffler where     $\Q'$ is the class of all pure submodules of modules in $\Q$.
\item[(ii)]
$M$ is
{\rm Prod}$\Q$-Mittag-Leffler.
\item[(iii)]
 $M$ is $\varinjlim\Q$-Mittag-Leffler.
\end{itemize}
\end{Th}
\begin{Proof}
(i) Let $\{Q_i\}_{i\in I}$ be a family of modules in $\Q$, and let
$\{Q'_i\}_{i\in I}$ be a family of left $R$-modules such that for
any $i\in I$ the module $Q'_i$ is a pure submodule of $Q_i$. For
any $i\in I$, denote by $\epsilon _i\colon Q'_i\to Q_i $ the
inclusion. As every $\epsilon _i$ is a pure monomorphism, so is
$\prod _{i\in I}\epsilon _i$. Then we have the commutative diagram
$$\begin{array}{ccccc}0&\to&M\otimes \prod _{i\in I}Q'_i & {\stackrel{M\otimes \prod _{i\in I}\epsilon _i
}{\longrightarrow}} & M\otimes \prod _{i\in I}Q_i\\
&&\rho '\downarrow \phantom{\rho '}& & \rho  \downarrow \phantom{\rho }\\
0&\to & \prod _{i\in I}(M\bigotimes_RQ'_i)&{\stackrel{\prod _{i\in
I}(M\otimes \epsilon _i)}{\longrightarrow}} & \prod _{i\in
I}(M\bigotimes_RQ_i).\end{array}
$$
As $\rho (M\otimes \prod _{i\in I}\epsilon _i)$ is injective, so
is $\rho '$.

(ii)
  is proved in \cite[p.~39]{rothmaler}.

(iii) We
follow the argument in \cite[Lemma~3.1]{facchini}.

 Let $\{Q_i\}_{i\in
I}$ be a family of modules such that, for any $i\in I$,
$Q_i=\varinjlim (K^i_\alpha, f^i_{\beta \, \alpha} )_{\beta \,
\alpha\in I  _i}$ and $K^i_\alpha\in \Q$ for any $\alpha \in I
_i$. For any $i\in I$ and $\alpha \in I   _i$, let $f^i_\alpha
\colon  K^i_\alpha \to Q_i$ denote the canonical morphism.

We want to show that $\rho \colon  M\bigotimes_R \prod _{i\in
I}Q_i\to \prod _{i\in I}(M\bigotimes_RQ_i)$ is injective.
Let $y=\sum _{j=1}^nx_j\otimes
(q^i_j)_{i\in I}$ be an element in the kernel of $\rho$.
This means that, for any $i\in I$, $\sum _{j=1}^nx_j\otimes q^i_j$
is the zero element of $M\otimes _RQ_i$. By Lemma~\ref{zeros}, for
each $i\in I$, there exists $\alpha _i\in I   _i$ and
$k^i_1,\dots ,k^i_n \in K^i_{\alpha_i}$ such that $\sum
_{j=1}^nx_j\otimes k^i_j$ is the zero element of $M\otimes
_RK^i_{\alpha_i}$ and $f^i_{\alpha_i}(k^i_j)=q^i_j$ for every
$j=1,\dots ,n$.

Consider the commutative diagram
$$\begin{array}{ccc}M\otimes \prod _{i\in I}K^i_{\alpha_i} & {\stackrel{M\otimes \prod _{i\in
I}f^i_{\alpha _i}
}{\longrightarrow}} & M\otimes \prod _{i\in I}Q_i\\
\rho '\downarrow \phantom{\rho '}& & \rho  \downarrow \phantom{\rho }\\
\prod _{i\in I}(M\bigotimes_RK^i_{\alpha_i})&{\stackrel{\prod
_{i\in I}(M\otimes f^i_{\alpha _i})}{\longrightarrow}} & \prod
_{i\in I}(M\bigotimes_RQ_i).\end{array}
$$
By construction,
\[y=\sum _{j=1}^nx_j\otimes
(q^i_j)_{i\in I}=(M\otimes \prod _{i\in I}f^i_{\alpha _i})(\sum
_{j=1}^nx_j\otimes (k^i_j)_{i\in I})\]
and
\[\rho '(\sum
_{j=1}^nx_j\otimes (k^i_j)_{i\in I})= (\sum _{j=1}^nx_j\otimes
k^i_j)_{i\in I}=0.\]
Note that $\rho'$ is injective because $K_{\alpha
_i }\in \Q$ for any $i\in I$. This shows that $y=0$.
\end{Proof}

\begin{Prop} \label{furtherclosures}
Let $R$ be a ring. The following statements hold true for    $M\in\Mod R$.

\begin{itemize}
\item[(i)] Let  $\Q_1, \ldots, \Q_n\subseteq R\LMod$, and let  $\Q =\cup _{i=1}^n\Q _i$.
If $M$ is $\Q_i$-Mittag-Leffler for all $1\le i\le n$, then $M$ is $\Q$-Mittag-Leffler.
\item[(ii)]  Let $\Q _1$ and $\Q _2$ be two classes in $R\LMod$, and let $\Q$ be the  class consisting
of all extensions of modules in $\Q _1$ by modules in $\Q _2$. Suppose that  $M$ is $\Q_i$-Mittag-Leffler for  $i=1,2$, and  that the functor $M\otimes -$ is exact on any short exact sequence with first term in $\Q_1$ and end-term in $\Q_2$.  Then $M$ is $\Q$-Mittag-Leffler.
\end{itemize}
\end{Prop}
\begin{Proof}
(i). Let $\{Q_i\}_{i\in I
}$ be a family of modules in $\Q$. For any $j=1,\dots, n$, set
\[I_j=\{i\in I\mid Q_i\in \Q_j \mbox{ and } Q_i\not\in \Q_k \mbox{ for } k<j\}.\]
Then $\prod _{i\in I}Q_i= \oplus _{j=1}^n(\prod _{i\in I_j} Q_i)$.
As $\rho_j\colon  M\otimes \prod _{i\in I_j} Q_i\to \prod _{i\in
I_j}(M\otimes  Q_i)$ is injective  for any $j=1,\dots n$, it
follows that $\rho$ is also injective.

(ii).    Let $\{Q_i\}_{i\in I}$ be a family
of left modules such that, for any $i\in I$, there is an exact
sequence
\[0\to Q_i^1\to Q_i\to Q_i^2\to 0\]
with $Q_i^1 \in \Q _1$ and $Q_i^2 \in \Q_2$. Then we have the
commutative diagram

$$\begin{array}{ccccccccc}&&M\otimes \prod _{i\in I}Q_i^1 & \to & M\otimes \prod _{i\in I}Q_i& {\to}&M\otimes \prod _{i\in I}Q_i^2&
\to&0\\ &&\rho _1\downarrow \phantom{\rho _1}& & \rho  \downarrow \phantom{\rho }&&\rho _2  \downarrow \phantom{\rho_2 }&&\\
0&\to & \prod _{i\in I}(M\bigotimes_RQ_i^1)&\to & \prod _{i\in
I}(M\bigotimes_RQ_i)&  {\to} & \prod _{i\in
I}(M\bigotimes_RQ_i^2)&\to&0\end{array}
$$
where the bottom row is exact by assumption on $M\otimes-$.
As $\rho _1$ and $\rho _2$ are injective,  $\rho$ is also
injective.
\end{Proof}

\begin{Cor}\label{tp}
Let $R$ be a ring, and  $M\in\Mod R$.
Let $({\mathcal T},{\mathcal F})$ be a torsion pair in $R\LMod$ such that   $M$ is ${\mathcal T}$-Mittag-Leffler and ${\mathcal F}$-Mittag-Leffler. Assume further that the functor $M\otimes -$ is exact on any short exact sequence with first term in ${\mathcal T}$ and end-term in ${\mathcal F}$.  Then $M$ is a Mittag-Leffler module.
\end{Cor}

\begin{Exs} \label{extp} {\rm
(1) Let $R$ be a commutative domain and denote by ${\mathcal T}$
and ${\mathcal F}$ the classes of torsion and torsionfree modules,
respectively.   Any flat $R$-module $M$ which is ${\mathcal
T}$-Mittag-Leffler and ${\mathcal F}$-Mittag-Leffler is a
Mittag-Leffler module.

(2) Let $\Lambda$ be a tame hereditary finite dimensional algebra over an algebraically closed field $k$, and let $\tube$ be the class of all finitely generated indecomposable regular modules. It was shown by Ringel in \cite[4.1]{Ri1}
that the classes $(\Fcal, \text{Gen}\tube)$ with $\Fcal=\tube^\circ={^\perp\tube}$ form a torsion pair, and for every module \mbox{$X\in \M$} there is a pure-exact sequence $$0\to \tube X\to X\to X/ \tube X\to 0$$ where
$\tube X=\sum_{f\in\Hom(Y,X), Y\in\tube}\im f \in \text{Gen}\,\tube$ is the trace of $\,\tube$ in $X$, and $X/ \tube X\in \Fcal$. Thus a module $M\in \M$ is Mittag-Leffler provided it is ${\mathcal T}$-Mittag-Leffler and ${\mathcal F}$-Mittag-Leffler.

(3) \cite[2.5]{Z} If $Q$ is a left $R$-module satisfying the
maximum   condition for finite matrix subgroups (see Definition
\ref{defmatrix}), for example an endonoetherian module, then every
right $R$-module is $Q$-Mittag-Leffler.

(4) \cite[2.4]{rothmaler}, \cite[2.1, 2.4]{Z} Let $\Q\subseteq R\LMod$. The class of
$\Q$-Mittag-Leffler modules is closed under pure submodules and pure extensions. A direct sum of
modules is $\Q$-Mittag-Leffler if and only if so are all direct summands. If $N$ is a finitely
generated  submodule of a $\Q$-Mittag-Leffler module $M$, then $M/N$ is $\Q$-Mittag-Leffler. }
\end{Exs}

Further examples are provided by the following results.

\begin{Prop} Let $R\to S$ be a ring epimorphism. Let $M_S$ be a finitely presented right
$S$-module, and let $N$ be a finitely generated $R$-submodule of $M$. Then $M_R$ and $M/N$ are
Mittag-Leffler with respect to the class $S$-$\mathrm{Mod}$.
\end{Prop}

\begin{Proof} Let $\{Q_i\}_{i\in I}$ be a family of left $S$-modules. Since $R\to S$ is a ring epimorphism
\[ M \otimes _R \prod _{i\in I} Q_i \cong M \otimes _R (S\otimes _S\prod
 _{i\in I} Q_i)= (M\otimes _R S)\otimes _S\prod
 _{i\in I} Q_i= M \otimes _S \prod _{i\in I} Q_i.\]
 As $M_S$ is finitely presented this is isomorphic to
\[\prod _{i\in I}M\otimes _S Q_i\cong \prod _{i\in I}M\otimes _RS\otimes _S
Q_i\cong \prod _{i\in I}M\otimes _R Q_i.\] This yields that the canonical map $\rho \colon M
\otimes _R \prod _{i\in I} Q_i \to \prod _{i\in I}M\otimes _R Q_i$ is in fact an isomorphism.

It follows from Example~\ref{extp}(4) that $M/N$ is also a Mittag-Leffler module with respect to
the class $S$-$\text{Mod}$. \end{Proof}

\begin{Def}\label{deffilters}
{\rm
Let $M$ be a right $R$-module, and  let $\tau$  denote an ordinal. An increasing chain $(M_\alpha\mid \alpha \le \tau)$
of submodules of $M$ is   a \emph{filtration} of $M$ provided that $M_0=0$, $M_\alpha =\bigcup
_{\beta <\alpha}M_\beta$ for all limit ordinals $\alpha \le \tau$ and $M_\tau =M$.

Given a class $\mathcal{C}$, a filtration $(M_\alpha\mid \alpha \le \tau)$ is a
\emph{$\mathcal{C}$-filtration} provided that $M_{\alpha +1}/M_\alpha \in \mathcal{C}$ for any
$\alpha < \tau$. We say also that $M$ is a \emph{$\mathcal{C}$-filtered module}.
}
\end{Def}

We have the following result about the  behavior of the Mittag-Leffler property with respect to
filtrations.

\begin{Prop}\label{MLfilters} Let $\mathcal{S}$ be a class of right $R$-modules that are Mittag-Leffler with respect
to a class $\mathcal{Q}\subseteq \mathcal{S}^\intercal$. Then any module isomorphic to direct
summand of an $\mathcal{S}\cup{\rm Add}\,R $-filtered module is $\Q$-Mittag-Leffler.
\end{Prop}

\begin{Proof} As projective modules are Mittag-Leffler and $(\mathcal{S}\cup{\rm Add}\,R )^\intercal= \mathcal{S}^\intercal$, we can assume that $\mathcal{S}$ contains ${\rm Add}\,R$ .
Moreover, since  the class of $\Q$-Mittag-Leffler modules is closed by direct
summands we only need to prove the statement for $\mathcal{S}$-filtered modules.

Let $M$ be an $\mathcal{S}$-filtered right $R$-module. Let $\tau$ be an ordinal such that there
exists an $\mathcal{S}$-filtration $(M_\alpha)_{\alpha \leq \tau}$ of $M$. We shall show that $M$
is $\Q$-Mittag-Leffler proving  by induction  that $M_\alpha $  is $\Q$-Mittag-Leffler for any
$\alpha \le \tau$. Observe that for any $\beta \leq \alpha \leq \tau$, $M_\alpha$ and $M_\alpha
/M_\beta$ are $\mathcal{S}$-filtered modules, so they belong to ${}^\intercal \mathcal{Q}$.

As $M_0=0$ the claim is true for $\alpha =0$. If  $\alpha <\tau$ then, as $\mathcal{Q}\subseteq
\mathcal{S}^\intercal$, we can apply
an argument similar to the one used in
Proposition~\ref{furtherclosures} to the exact sequence
\[0\to M_\alpha \to M_{\alpha +1}\to M_{\alpha +1}/M_\alpha \to 0\]
to conclude that if $M_\alpha$ is $\Q$-Mittag-Leffler then $M_{\alpha +1}$ is $\Q$-Mittag-Leffler.

Let $\alpha \leq \tau$ be a limit ordinal, and assume that $M_\beta$ is $\Q$-Mittag-Leffler for
any $\beta <\alpha$. We shall prove that $M_\alpha =\bigcup _{\beta <\alpha}M_\beta$ is
$\Q$-Mittag-Leffler. Let $\{Q_i\}_{i\in I}$ be a family of modules in $\Q$, and let $x\in
\mathrm{Ker}\, (M_\alpha \otimes _R\prod _{i\in I}Q_i\to \prod _{i\in I} M_\alpha \otimes _R
Q_i)$. There exists $\beta <\alpha $ and $y\in M_\beta \otimes _R \prod _{i\in I}Q_i$ such that
$x=(\epsilon _\beta \otimes _R \prod _{i\in I}Q_i)\, (y)$, where $\epsilon _\beta \colon M_\beta \to
M_\alpha$ denotes the canonical inclusion. Considering the commutative diagram
$$\begin{array}{ccccc}&&M_\beta \otimes \prod _{i\in I}Q_i & {\stackrel{\epsilon_\beta \otimes \prod _{i\in
I}Q_i}
{\longrightarrow}} & M_\alpha \otimes \prod _{i\in I}Q_i\\
&&\rho'\downarrow \phantom{\rho'}& & \rho  \downarrow \phantom{\rho }\\
&& \prod _{i\in I} (M_\beta  \otimes  Q_i)&{\stackrel{\prod _{i\in I} (\epsilon _\beta \otimes
Q_i)}{\longrightarrow}} & \prod _{i\in I} (M_\alpha \otimes Q_i)\end{array}
$$
we see that $0=\prod _{i\in I} (\epsilon _\beta \otimes Q_i)\rho '(y)$. As $\rho '$ is injective
because $M_\beta$ is $\Q$-Mittag-Leffler and, for any $i\in I$, $\epsilon _\beta \otimes Q_i$ is
also injective because $\mathrm{Tor}^R_1(M_\alpha /M_\beta, Q_i)=0$, we deduce that $y=0$.
Therefore $x=0$, and $\rho$ is injective.
\end{Proof}

\section{$H$-subgroups}\label{H}

We  recall a notion from \cite{zbcn} which will be very useful in the sequel.

\begin{Def} {\rm Let $M$, $M'$ and $B$ be right $R$-modules, and let $v\in \Hom
_R(M,M')$. The subgroup (and End$_RB$-submodule) of $\Hom _R(M,B)$ consisting of all compositions of $v$ with maps
in $\Hom _R(M',B)$ is denoted by  $$H_v(B) =\Hom _R(M',B)\,v$$ and is called an
\emph{$H$-subgroup} of $\Hom _R(M,B)$.
}
\end{Def}

\begin{Remark}\label{hsubgr} {\rm  \cite[2.10]{zbcn}
Let $M$, $M'$ and $B$ be right $R$-modules, and let $v\in \Hom
_R(M,M')$. \\
(1) $H_v$ is a subfunctor of $\Hom _R(M,-)$ commuting with direct products. If $M$ is finitely generated, then  $H_v$ also commutes with direct sums.\\
(2)  An $\End B$-submodule $U$ of $\Hom _R(M,B)$ is an $H$-subgroup if and only if there are a set $I$ and a homomorphism $u\in \Hom
_R(M,B^I)$ such that $U=H_u(B)$.
}
\end{Remark}

\noindent
For the following discussion it is important to keep in mind the
following easy observations.

\begin{Lemma}\label{diagram}
Let $M$, $M'$, $N$  be right $R$-modules,   $u\in\Hom _R(M,N)$,   $v\in \Hom
_R(M,M')$.
\begin{enumerate}
\item
If there is  $h\in \Hom _R(M',N)$ such that the diagram
\[\begin{array}{ccc}
M & {\stackrel{v}{\longrightarrow}}&M'\\
{\scriptstyle u} \downarrow \phantom{u }&{\phantom{h}{\swarrow}}{{}_h} & \\
N&&
\end{array}\]
 commutes, then
$H_u(B)\subseteq H_v(B)$ for any right $R$-module $B$.
\item  If $B$ is a right $R$-module such that $H_u(B)\subseteq H_v(B)$, then $H_{ut}(B)\subseteq H_{vt}(B)$ for all $t\in\Hom _R(X,M)$, $X\in \Mod R$.
\end{enumerate}
\end{Lemma}

\noindent
Recall that a homomorphism $\pi:B\to B''$ is a
{\it locally split epimorphism}
if for each finite subset $X\subseteq B''$ there is a map $\varphi=\varphi_X:B''\to B$ such that $x=\pi\varphi(x)$ for all $x\in X$.
Observe that  every split epimorphism  is locally split, and every locally split epimorphism  is a pure epimorphism.
{\it Locally split}   {\it monomorphisms}   are defined dually.
Moreover, a submodule $B'$ of a module $B$ is said to be a
{\it locally split} (or strongly pure \cite{Zlpi}) {\it  submodule}
if the embedding $B'\subset B$ is locally split.

\begin{Lemma}\label{star}
Let $M$ and $M'$  be right $R$-modules, and let $v\in \Hom _R(M,M')$. Assume that $M$ is finitely
generated. Let further $\varepsilon: B'\to B$ be a pure monomorphism. If $M'$ is finitely
presented or $\varepsilon$ is a locally split monomorphism, then \[\varepsilon H_v(B')=H_v(B)\cap
\varepsilon \mathrm{Hom}_R(M,B').\]
\end{Lemma}
\begin{Proof}
The first case is treated in \cite[Lemma~4.1]{bazher} or \cite[Lemma~2.8]{abh}. For the second
case, we assume that $\varepsilon$ is a locally split monomorphism. We show the inclusion
$\supseteq$.  Consider $f\in \mathrm{Hom}_R(M,B')$ such that  $\varepsilon\,f=h\,v$ for some $h\in
\mathrm{Hom}_R(M',B)$. Choose a generating set $x_1,\ldots,x_n$    of $M$ together with a map
$\varphi:B\to B'$ such that $f(x_i)=\varphi\varepsilon\,f(x_i)$ for all $1\le i\le n$. Then the
composition $h'=\varphi\, h$ satisfies $f=h'\,v\in H_v(B')$. The inclusion $\subseteq$ is clear.
\end{Proof}

\begin{Lemma}\label{hpure} Assume that  the diagram
\[\begin{array}{ccc}
M & {\stackrel{v}{\longrightarrow}}&M'\\
{\scriptstyle u}  \downarrow \phantom{u }&{\phantom{h}{\swarrow}}{{}_h} & \\
N&&
\end{array}\]
of  right $R$-modules and module homomorphisms
commutes. Assume further that $B$ is a right $R$-module such that $H_u(B)= H_v(B)$.
\begin{enumerate}
\item
If $h$ factors through a homomorphism $m\in\Hom _R(M',M'')$, then $H_u(B)= H_{mv}(B)$.
\item Assume that $M$ is finitely generated and $N$ is finitely presented.
If $B'\subseteq B$ is a pure submodule, then $H_u(B')= H_v(B')$.
\item Assume that $M$ is finitely generated and $M'$ is finitely presented.
If $ B
\stackrel{\pi}\to B''$ is a pure-epimorphism , then
$H_u(B'')= H_v(B'')$.
\item If $M$ is finitely generated and $B'\subseteq B$ is a locally split  submodule, then $H_u(B')= H_v(B')$.
\item If $M$ is finitely generated and $ B
\stackrel{\pi}\to B''$ is a locally split epimorphism,  then
$H_u(B'')= H_v(B'')$.
\end{enumerate}
\end{Lemma}

\begin{Proof}
(1) is left to the reader. For the remaining statements, note first
that by Lemma~\ref{diagram} it suffices to show  $H_v(B')\subseteq H_u(B')$ and $H_v(B'')\subseteq H_u(B'')$, respectively.
For
(2) and (4), observe that
 Lemma \ref{star} yields
$\varepsilon H_u(B')=H_u(B)\cap \varepsilon \mathrm{Hom}_R(M,B') $ where $\varepsilon:B'\to B$ denotes the canonical embedding. Then
$\varepsilon H_v(B')\subseteq H_v(B)\cap \varepsilon
\mathrm{Hom}_R(M,B')=
\varepsilon H_u(B')$, and since $\varepsilon$ is a
monomorphism, we deduce that $H_v(B')\subseteq H_u(B')$.
\\
In statement (3), we have that
  $\pi \colon B\to B''$ is a pure epimorphism and $M'$ is
finitely presented, so \[\mathrm{Hom}_R(M',\pi)\colon
\mathrm{Hom}_R(M',B)\to \mathrm{Hom}_R(M',B'')\] is also an
epimorphism. Thus, if $fv\in H_v(B'')$, then there exists $g\in
\mathrm{Hom}_R(M',B)$ such that $\pi g=f$. By hypothesis
$gv\in H_u(B)$, so $fv=\pi g v \in \pi H_u(B)\subseteq H_u(B'')$.\\
For statement (5), we consider $fv\in H_v(B'')$, and choose a generating set $x_1,\ldots,x_n$    of $M$ together with a map
$\varphi:B''\to B$ such that $fv(x_i)=\pi\varphi\,fv(x_i)$ for all $1\le i\le n$.
Then the composition $h=\varphi\, f$ satisfies $fv=\pi\,hv$. Moreover, $h\,v\in H_u(B)$, so there is $h'\in \Hom _R(N,B)$ such that $hv=h'u$. Thus
$fv=\pi\,h'u\in H_u(B'').$
\end{Proof}


\section{$\B$-stationary modules}
\label{stat}

\begin{Def} {\rm An inverse system of sets $(H_\alpha,h_{\alpha\,\gamma})_{\alpha\,\gamma\in
I  }$ is said to satisfy the \emph{Mittag-Leffler condition} if for
any $\alpha \in I  $ there exists $\beta\ge \alpha$
such that
$h_{\alpha\,\gamma}(H_\gamma)= h_{\alpha\,\beta}(H_{\beta})$  for any $\gamma \ge \beta$.}
\end{Def}

Let us specify the Mittag-Leffler condition for the case $I   =\N$.

\begin{Ex} \label{countableML} An inverse system of the form
\[\cdots H_{n+1}\stackrel{h_n}\to H_n\cdots H_{2}\stackrel{h_1}\to
H_1\] satisfies the Mittag-Leffler condition if and only if for
any $n\in \N$ the chain of subsets of $H_n$
\[h_n(H_{n+1})\supseteq \cdots \supseteq h_n\cdots h_{n+k}(H_{n+k+1})\supseteq \cdots\]
is stationary.
\end{Ex}

According to Raynaud and Gruson \cite[p.~74]{RG} the following
characterization of Mittag-Leffler inverse systems is due to
Grothendieck as it is implicit in \cite[13.2.2]{grothendieck}. We
give a proof for completeness' sake.

\begin{Lemma}\label{MLstrict} Consider an inverse system of the form
\[\mathcal{H} \colon \cdots H_{n+1}\stackrel{h_n}\to H_n\cdots H_{2}\stackrel{h_1}\to
H_1.\] For any $m>n\ge 1$ set $h_{nm}=h_n\cdots h_{m-1}$, and, for
any $n\ge 1$ let $g_n\colon \varprojlim H_i\to H_n$ denote the
canonical map.

The inverse system $\mathcal{H}$ satisfies the Mittag-Leffler
condition if and only if for any $n\ge 1$ there exists $\ell (n)>n$
such that
\[g_n(\varprojlim H_i)=h_{n\, \ell (n)}(H_{\ell (n)})= h_n\cdots
h_{\ell(n)-1} (H_{\ell (n)}).\]
\end{Lemma}

\begin{Proof} Observe that since, for any $m>n \ge 1$, $g_n=h_{n
m}g_m$ always
\[g_n(\varprojlim H_i)\subseteq \bigcap _{m>n} h_{nm}(H_{m}).\]

Assume now that $\mathcal{H}$ satisfies the Mittag-Leffler
condition. We only need to show that for any $n\ge 1$ there exists
$\ell (n)> n$ such that $h_{n\, \ell (n)}(H_{\ell (n)})\subseteq g_n(\varprojlim
H_i).$ To this aim fix $n\ge 1$.

Applying repeatedly that $\mathcal{H}$ satisfies the Mittag-Leffler
condition we find a sequence of elements in $\mathbb{N}$
\[(*) \qquad n=n_0<n_1<\cdots <n_i< \cdots\]
such that $h_{n_i\,
n_{i+1}}(H_{n_{i+1}})=h_{n_i\, m}(H_m)$ for all $i\ge 0$ and $m\ge n_{i+1}$. Now we show that $\ell (n)$
can be taken to be $n_1$.

Let $a\in h_{n_0 \, n_1}(H_{n_1})$. Then $a\in h_{n_0 n_2}(H_{n_2})$, and there is
$a_1\in  h_{n_1 n_2}(H_{n_2})\subseteq H_{n_1}$ such that $a=h_{n_0 n_1}(a_1)$. In this fashion, the
 properties of the sequence
$(*)$ allow us to find a sequence $a_0=a, a_1,\dots ,a_i ,\dots $
such that $a_i\in H_{n_i}$ and $h_{n_i\,n_{i+1}}(a_{i+1})=a_i$ for
any $i\ge 0$. Hence $b=(a_i) \in \varprojlim
H_{n_i}=\varprojlim H_j$ and $g_n(b)=a_0=a$ as desired.

The converse implication is clear because of the remarks at the
beginning of the proof.
\end{Proof}

The above characterization does not extend to uncountable inverse
limits; an example where this fails is implicit in
Example~\ref{notstrict}.

\medskip

 We will be interested in inverse systems arising by applying the functor $\Hom _R(-,B)$ on a direct system.

\begin{Remark}\label{translation} Let $(F_\alpha,u_{\beta \,\alpha})_{\beta \,\alpha \in I  }$
be a direct system of right $R$-modules,  $B$   a right
$R$-module, and $\beta \ge \alpha \in I  $. Then  $$(\Hom _R(F_\alpha, B), \Hom _R(u_{\beta
\,\alpha}, B))_{\beta \,\alpha \in I  }$$ is an   inverse system of left modules over the endomorphism ring of $B$, and
  $$\Hom _R(u_{\beta
\,\alpha},B)(\Hom _R(F_{\beta},B))=\Hom _R(F_{\beta},B)\,u_{\beta
\,\alpha}=H_{u_{\beta \,\alpha}}(B)$$
\end{Remark}

Applying Lemma~\ref{diagram}(1) to the situation of Remark \ref{translation} we obtain the following.

\begin{Lemma}\label{Hinversesystem}  Let $(F_\alpha,u_{\beta \,\alpha})_{\beta \,\alpha \in
I  }$ be a direct system of right $R$-modules  with direct limit
$M$, and denote by $u_\alpha \colon  F_\alpha \to M$ the canonical
map. Let $B$ be a right $R$-module.
\begin{itemize}
\item[\textrm{(i)}] If $\gamma \ge \beta \ge \alpha$ then $H_{u_{\gamma \, \alpha}}(B)\subseteq H_{u_{\beta \,
\alpha}}(B)$.
\item[\textrm{(ii)}] $H_{u_{\alpha}}(B)\subseteq H_{u_{\beta \,
\alpha}}(B)$ for any $\beta \ge \alpha$.
\end{itemize}
\end{Lemma}

This allows to interpret the Mittag-Leffler condition on inverse systems as in Remark \ref{translation} in terms of $H$-subgroups.
\begin{Lemma}\label{alphabeta} Let $(F_\alpha,u_{\beta \,\alpha})_{\beta \,\alpha \in
I  }$ be a direct system of right $R$-modules. Let $\alpha$,
$\beta \in I  $ with $\beta\ge\alpha$, and let $B$ be  a right $R$-module. The
following statements are equivalent.
\begin{itemize}
\item[(1)] For any $\gamma\ge \alpha$, the inclusion $H_{u_{\beta\, \alpha}}(B)\supseteq H_{u_{\gamma\,
\alpha}}(B)$ implies $H_{u_{\beta\, \alpha}}(B)= H_{u_{\gamma\,
\alpha}}(B)$.
\item[(2)] $H_{u_{\beta\, \alpha}}(B)\subseteq \bigcap _{\gamma \ge \beta} H_{u_{\gamma\,
\alpha}}(B)$.
\item[(2')] $H_{u_{\beta\, \alpha}}(B)= \bigcap _{\gamma \ge \beta} H_{u_{\gamma\,
\alpha}}(B)$.
\item[(3)] $H_{u_{\beta\, \alpha}}(B)\subseteq \bigcap _{\gamma \ge \alpha} H_{u_{\gamma\,
\alpha}}(B)$.
\item[(3')] $H_{u_{\beta\, \alpha}}(B)= \bigcap _{\gamma \ge \alpha} H_{u_{\gamma\,
\alpha}}(B)$.
\end{itemize}
\end{Lemma}
\begin{Proof} By Lemma~\ref{Hinversesystem} (i) it follows
immediately that $(1)$ implies $(2)$, and that $(2)$ and $(2')$,
as well as $(3)$ and $(3')$ are equivalent statements. Further, it is clear that $(3)\Rightarrow (1)$.

We show $(2)\Rightarrow (3)$. Let $\gamma \ge \alpha$ and choose
$\gamma _1\in I  $ such that $\gamma _1\ge \gamma$ and $\gamma
_1\ge \beta$. By $(2)$, $H_{u_{\beta\, \alpha}}(B)\subseteq
H_{u_{\gamma _1\, \alpha}}(B)$ and, by Lemma~\ref{Hinversesystem}
(i), $H_{u_{\gamma _1\, \alpha}}(B)\subseteq H_{u_{\gamma \,
\alpha}}(B)$. Hence, $H_{u_{\beta\, \alpha}}(B)\subseteq \bigcap
_{\gamma \ge \alpha} H_{u_{\gamma\, \alpha}}(B)$ as we wanted to
proof.
\end{Proof}

We adopt the following definition inspired by the terminology in \cite{GH}.

\begin{Def}\label{Bstat}{\rm Let $B$ be a  right
$R$-module.

(1) A direct system $(F_\alpha,u_{\beta \,\alpha})_{\beta \,\alpha \in I  }$
 of right $R$-modules  is said to be \emph{$B$-stationary} provided that the inverse system $(\Hom _R(F_\alpha, B), \Hom _R(u_{\beta
\,\alpha}, B))_{\beta \,\alpha \in I  }$
satisfies the Mittag-Leffler condition, in other words, provided for any $\alpha \in I  $ there exists $\beta \ge
\alpha$ such that  the equivalent
conditions in Lemma~\ref{alphabeta} are satisfied.

(2) A right $R$-module $M$ is said to be  \emph{$B$-stationary} if there exists a {$B$-stationary} direct system of
finitely presented modules $(F_\alpha,u_{\beta \,\alpha})_{\beta \,\alpha \in I  }$ such that  $M=\varinjlim F_\alpha$.

(3) Let $\B$ be a class of right
$R$-modules. We say that a direct system $(F_\alpha,u_{\beta \,\alpha})_{\beta \,\alpha \in I  }$
or a right $R$-module $M$ are $\B$-stationary if they are $B$-stationary for all $B\in\B$.}
\end{Def}

Let us start by discussing some closure properties of the class $\B$.

\begin{Prop}\label{prodstationary} Let $\{B_j\}_{j\in J}$ be a
family of right $R$-modules. Let $(F_\alpha,u_{\beta
\,\alpha})_{\beta \,\alpha \in I  }$ be a direct system of right
$R$-modules. Then the following statements are equivalent.
\begin{itemize}
\item[(1)] $(F_\alpha,u_{\beta
\,\alpha})_{\beta \,\alpha \in I  }$ is $\prod _{j\in
J}B_j$-stationary.
\item[(2)] For any $\alpha \in I$ there exists $\beta \ge \alpha$ such
that $H_{u_{\beta \alpha}}(B_j)=\bigcap _{\gamma \ge
\alpha}H_{u_{\gamma \alpha}}(B_j)$ for any $j\in J$.
\end{itemize}
If $F_\alpha$ is finitely generated for any $\alpha \in I$, then
the above statements are further equivalent to
\begin{itemize}
\item[(3)] $(F_\alpha,u_{\beta
\,\alpha})_{\beta \,\alpha \in I  }$ is $\bigoplus _{j\in
J}B_j$-stationary.
\end{itemize}
\end{Prop}

\begin{Proof} We use the same arguments as in \cite[2.6]{abh}.
$(1)\Leftrightarrow (2).$ Statement $(1)$ holds if and only if, for any $\alpha \in I$
there exists $\beta$ such that $H_{u_{\beta \alpha}}(\prod _{j\in
J}B_j)=\bigcap _{\gamma \ge \alpha}H_{u_{\gamma \alpha}}(\prod _{j\in
J}B_j)$ if and only if
\[  \prod _{j\in
J}H_{u_{\beta \alpha}}(B_j)=\bigcap _{\gamma \ge \alpha}\prod _{j\in
J} H_{u_{\gamma \alpha}}(B_j)= \prod _{j\in J} \bigcap _{\gamma \ge
\alpha}H_{u_{\gamma \alpha}}(B_j).\] Equivalently, if and only if
$(2)$ holds.

The proof of $(2)\Leftrightarrow (3)$ follows in a similar way by
observing that
\[H_{u_{\beta \alpha}}(\bigoplus _{j\in J}B_j)=\bigoplus _{j\in J}H_{u_{\beta \alpha}}(B_j).\]
provided all $F_\alpha$ are finitely generated.
\end{Proof}

\begin{Cor} \label{closureB} Let $\B$ be a class
of right $R$-modules. Let $M$ be a $\B$-stationary  right
$R$-module. Then the following statements hold true.
\begin{itemize}
\item[(i)] $M$ is $\B '$-stationary where $\B '$ denotes the class
of all modules isomorphic either to a pure submodule  or to a
pure quotient  of a module  in $\B$.

\item[(ii)] $M$ is $\Add \,\B$-stationary if and only if it
is $\Prod\, \B$-stationary if and only if there exists a direct
system of finitely presented right $R$-modules
$(F_{\alpha},u_{\beta \alpha})_{\beta \alpha \in I}$ with
$\varinjlim F_{\alpha}\cong M$ having the property that for any $\alpha \in I$  there
exists $\beta \ge \alpha$ such that $H_{u_{\beta \alpha}}(B)=\bigcap
_{\gamma \ge \alpha}H_{u_{\gamma \alpha}}(B)$ for any $B\in \B$.

\item[(iii)] $M$ is $\Add B$- and $\Prod B$-stationary for every $B\in \B$.

\end{itemize}
\end{Cor}

\begin{Proof}
 The statements in Lemma~\ref{hpure}(2) and (3) imply statement $(i)$. Statement
 $(ii)$ is a direct consequence of   Proposition~\ref{prodstationary}  combined with $(i)$, and $(iii)$ is a special case of $(ii)$. \end{Proof}





\begin{Prop}\label{HML} Let $(F_\alpha,u_{\beta \,\alpha})_{\beta \,\alpha \in
I  }$ be a direct system of right $R$-modules, and let $B$ be  a right $R$-module. Consider
the following  statements.
\begin{itemize}
\item[(1)]
For any infinite chain $\alpha _1\le \alpha _2\le \dots \in I  $ the  direct system $(F_{\alpha_n},u_{\alpha _{n+1} \,\alpha _n})_{n \in \N}$ is $B$-stationary.

\item[(1')] For any infinite chain $\alpha _1\le \alpha _2\le \dots \in I  $
the chain of subgroups
\[H_{u_{\alpha _2\,\alpha _1}}(B)\supseteq H_{u_{\alpha _3\,\alpha _1}}(B)\supseteq \cdots\]
is stationary.

\item[(2)] The direct system  $(F_\alpha,u_{\beta \,\alpha})_{\beta \,\alpha \in
I  }$ is $B$-stationary.

\end{itemize}
Then $(1)$ and $(1')$ are equivalent statements which imply $(2)$.
\end{Prop}
\begin{Proof} The fact that $(1)$ and $(1')$
are equivalent statements follows directly from the definitions
taking into account Example~\ref{countableML} and
Remark~\ref{translation}.

We prove now $(1)\Rightarrow (2)$. Assume for a
contradiction that there exists $\alpha$ such that for any
$\beta\ge \alpha$ condition $(1)$ in Lemma~\ref{alphabeta} fails.
Now we construct a countable chain in $I  $ such that
condition $(1)$ fails.

Set $\alpha _1=\alpha$. Let $n\ge 1$, and assume we have
constructed $\alpha _1\le \alpha _2 \dots \le \alpha _n$ such that
\[H_{u_{\alpha _2\,\alpha _1}}(B)\varsupsetneq H_{u_{\alpha _3\,\alpha _1}}(B)\cdots \varsupsetneq
H_{u_{\alpha _n\,\alpha _1}}(B).\]

As Lemma~\ref{alphabeta} $(1)$
fails for $\alpha _n\ge \alpha _1$, there exists $\gamma \ge
\alpha _1$ such that
\[H_{u_{\alpha _n\,\alpha _1}}(B)\varsupsetneq H_{u_{\gamma\,\alpha _1}}(B).\]
Let $\alpha _{n+1}\in I  $ be such that $\alpha _{n+1}\ge
\gamma$ and $\alpha _{n+1}\ge \alpha _n$. By
Lemma~\ref{Hinversesystem} (i),
\[H_{u_{\alpha _{n+1}\,\alpha _1}}(B)\subseteq H_{u_{\gamma\,\alpha _1}}(B)\subsetneq H_{u_{\alpha _n\,\alpha _1}}(B)\]
as wanted.
\end{Proof}

For later reference, we recall the following result.

\begin{Th}{\rm \cite[Theorem~5.1]{bazher}} \label{BH}
Let $\B$ be a class of right $R$-modules such that if $B\in \B$
then $B^{(\N)}\in \B$, and let $\A={}^\perp \B$. Let moreover
\[F_1\stackrel{u_1}\to F_2\stackrel{u_2}\to F_3\to\dots\to F_{n}\stackrel{u_n}\to F_{n+1}\to\dots\]
be a countable direct system of  finitely presented right
$R$-modules, and consider the   pure exact sequence
\[(\ast)\hspace{1cm} 0\to \oplus_{n\in \N} F_n\stackrel{\phi}\to \oplus_{n\in \N} F_n\to \displaystyle \lim_{\longrightarrow} F_n\to0\]
where $\phi\varepsilon_n=\varepsilon_n -\varepsilon_{n+1}u_n$ and
$\varepsilon_n\colon  F_n\to \oplus_{n\in \N} F_n$ denotes the
canonical morphism for every $n \in \N$. Then the following
statements are equivalent.
\begin{itemize}
\item[(1)] The direct system $(F_n,u_n)_{n\in\N}$ is $\B$-stationary.
\item[(2)] $\mathrm{Hom}_R(\phi, B)$ is surjective for all $B\in\B$.
\item[(3)] $\varprojlim^1 \Hom _R(F_n,B)=0$ for all $B\in\B$.
\end{itemize}
If $F_n$ belongs to  $\A$ for all $n\in\N$, then the following statement is further equivalent.
\begin{itemize}
\item[(4)] $\varinjlim F_n \in\A$.
\end{itemize}
\end{Th}

\medskip

\begin{Cor}\label{reducetocountable} Let $\B$ be a class of right $R$-modules such that if $B\in \B$ then
$B^{(\N)}\in \B$, and let $\A={}^\perp \B$ .Then the following
statements are equivalent.
\begin{itemize}
\item[(1)] Every  countable direct system of finitely presented
modules in $\A$ has limit in $\A$.
\item[(2)] Every countable direct system of finitely presented modules in $\A$ is $\B$-stationary.
\item[(3)] Every direct system of finitely presented modules in $\A$ is $\B$-stationary.
\end{itemize}
\end{Cor}

\begin{Proof} Assume $(1)$. Let $(F_\alpha,u_{\beta \,\alpha})_{\beta \,\alpha \in
I  }$ be a direct system of finitely presented right $R$-modules
such that $F_\alpha \in \A $ for any $\alpha \in I  $. Let $\alpha
_1\le \alpha _2\le \dots $ be a chain in $I  $. By $(1)$,
$\varinjlim (F_{\alpha _n},u_{\alpha_{n+1}\, \alpha _n})\in \A$.
Then $(F_{\alpha _n},u_{\alpha_{n+1}\, \alpha _n})_{n\in \N}$  is
$\B$-stationary by \ref{BH}, hence condition $(3)$ follows by
Proposition~\ref{HML}.

Obviously $(3)$ implies $(2)$. To see that $(2)$ implies $(1)$,
let $A=\varinjlim (F_\alpha,u_{\beta \,\alpha})_{\beta \,\alpha
\in I }$ be such that $I  $ is countable and $F_\alpha$ are
finitely presented modules in $\A$. Taking a cofinal set of $I  $
if necessary we may assume that $I   =\N$. Our hypothesis allows
us to use \ref{BH} to conclude that $A\in \A$.
\end{Proof}

\begin{Exs} \label{exstat} {\rm
(1) Let $\B$ be a class of right $R$-modules such that if $B\in
\B$ then $B^{(\N)}\in \B$, and let $\A= {}^\perp \B$. If $M\in \A$
is countably presented,  then $M$ is $\B$-stationary.

In fact, $M$ can be written as direct limit of a countable direct system as in \ref{BH}, and for all modules $B\in\B$ the map $\mathrm{Hom}_R(\phi, B)$ is surjective  because
$\Ext _R^1 (\varinjlim F_n, B)=0$.

\medskip

(2) Let $\B$ and $\A$ be as in (1). Assume that
 $R$ is a right  noetherian
ring and $\B$ consists of modules of injective dimension at most one. Then every $M\in \A$
is $\B$-stationary.

In fact, the additional assumption on $\B$ means that $\A$ is closed by submodules: Let $N\leq M\in
\A$. For any $B\in \B$, if we apply $\Hom _R(-,B)$ to the exact
sequence
\[0\to N\to M\to M/N\to 0\]
we obtain the exact sequence \[\Ext _R^1(M/N,B)\to \Ext
_R^1(M,B)=0\to \Ext _R^1(N,B)\to \Ext _R^2(M/N,B)=0 .\] Hence,
$\Ext _R^1(N,B)=0$.

As $R_R$ is noetherian, any finitely generated submodule of $M$ is
finitely presented. Let $I  $ denote the directed set of all
finitely generated submodule of $M$, then $M=\bigcup _{F\in
I  }F$. If $F_1\le F_2\le \dots \le F_n\le \dots$ is a chain
in $I  $, then $N=\bigcup _{n\in \N}F_n$ is a submodule of $M$
and it is in $\A$. By \ref{BH}, $N$ is
$\B$-stationary. By Proposition~\ref{HML}, $M$ is $\B$-stationary.

\medskip

(3)  Let $M$ be a module with a perfect decomposition in the sense of \cite{Angsaor},
for example $M$ a $\Sigma$-pure-injective module, or $M$ a finitely generated module with perfect endomorphism ring.
Let $\eme$ be a class of finitely presented modules in Add$M$. Then every $N\in\varinjlim\eme$ is $\Mod R$-stationary.

In fact,  we can write $N=\varinjlim F_\alpha$ where
$(F_\alpha,u_{\beta \,\alpha})_{\beta \,\alpha \in I  }$ is a direct
system of finitely presented modules in $\eme$. If we take a chain
$\alpha _1\le \alpha _2\le \dots $  in $I  $, then  $(F_{\alpha
_n},u_{\alpha_{n+1}\, \alpha _n})_{n\in \N}$ is a direct system in
Add$M$ with a totally ordered index set, so it follows from
\cite[1.4]{Angsaor} that the pure exact sequence  $(\ast)$
considered in Theorem \ref{BH} is split exact. In particular,
$\mathrm{Hom}_R(\phi, B)$ is surjective for all modules $B$, hence
$(F_{\alpha _n},u_{\alpha_{n+1}\, \alpha _n})_{n\in \N}$  is $\Mod
R$-stationary by \ref{BH}. Now the claim  follows from
Proposition~\ref{HML}.

\medskip

(4) Let $B$ be a $\Sigma$-pure-injective module. Then every right $R$-module $M$ is Add$B$-stationary.

To see this, write $M=\varinjlim F_\alpha$ where $(F_\alpha,u_{\beta
\,\alpha})_{\beta \,\alpha \in I  }$ is a direct system of finitely
presented modules. If we take a chain $\alpha _1\le \alpha _2\le
\dots $ in $I  $ and consider the direct system  $(F_{\alpha
_n},u_{\alpha_{n+1}\, \alpha _n})_{n\in \N}$, then for any
$B'\in\text{Add}B$ we know that  $\mathrm{Hom}_R(-, B')$ is exact on
the pure exact sequence  $(\ast)$ considered in Theorem \ref{BH}. So
$\mathrm{Hom}_R(\phi, B')$ is surjective for all modules
$B'\in\text{Add}B$, and the claim follows again by combining Theorem
\ref{BH} and Proposition~\ref{HML}. }
\end{Exs}

\section{Dominating maps}\label{dom}

{}From the characterization of Mittag-Leffler modules in \cite{RG}, we know  that a right module
is
 $Q$-Mittag-Leffler for any left module $Q$ if and only if
 it is $B$-stationary for any right module $B$. We will now  investigate the relationship between
the properties $Q$-Mittag-Leffler  and $B$-stationary when we restrict our choice of $Q$ and $B$ to  subclasses
  of $R\LMod$ and $\Mod R$, respectively.

As a first step, in Theorem \ref{2.1.4} we provide a characterization of when a module $M$ is  $B$-stationary
which is independent from the direct limit presentation of $M$. To this end, we need the following notion
which is inspired by the corresponding notion from \cite{RG}.

\begin{Def}\label{dominate} {\rm Let $Q$ be a left $R$-module, and let $B$ be a right $R$-module.
Let moreover $u\colon  M\to N$ and $v\colon  M\to M'$ be right
$R$-module homomorphisms.  We say that $v$ \emph{$B$-dominates $u$
with respect to $Q$} if
\[\Ker (u\otimes _RQ)\subseteq \bigcap _{h\in H_v(B)} \Ker (h\otimes _RQ) \]
For classes of modules $\Q$ and $\B$ in $R\LMod$ and $\Mod R$, respectively,  we   say that $v$ \emph{$\B$-dominates $u$
with respect to $\Q$} if $v$ $B$-dominates $u$ with respect to $Q$
 for any $Q\in\Q$ and any $B\in \B$.

 If $\Q=R$-$ \mathrm{Mod}$, we simply say that $v$ $\B$-dominates
 $u$.
 If $\B=\mathrm{Mod}$-$R$, we say that $v$  dominates $u$ with respect to
 $\Q$, and
 of course, this means that  $\Ker (u\otimes _RQ)\subseteq \Ker (v\otimes _RQ)$ for all left modules $Q\in \Q$.

 Finally, if $\Q=R$-$ \mathrm{Mod}$ and $\B=\mathrm{Mod}$-$R$, then we are in the case treated  in \cite[2.1.1]{RG}, and we say that $v$  dominates
 $u$.
}
\end{Def}

We note some properties of  dominating maps.

\begin{Lemma}\label{composition} Let $u\colon  M\to N$ and $v\colon  M\to M'$ be right
$R$-module homomorphisms, and let $B$ be a right $R$-module and
$Q$ a left $R$-module.

\begin{itemize}
\item[(1)] ``$B$-dominating with respect to $Q$'' is translation invariant on the right. That
is, if $v$ $B$-dominates $u$ with respect to $Q$ and $t\colon X\to
M$ is a homomorphism, then $vt$ $B$-dominates $ut$ with respect to
$Q$.

\item[(2)] ``$B$-dominating with respect to $Q$'' is stable by composition on the left. More
precisely, if $v$ $B$-dominates $u$ with respect to $Q$ and
$m\colon  M'\to M''$ is a homomorphism, then $mv$ $B$-dominates
$u$ with respect to $Q$.
\end{itemize}
\end{Lemma}
\begin{Proof} (1) By hypothesis, $\Ker (ut\otimes Q)=\Ker (u\otimes Q)(t\otimes
Q)$ is contained in
\[\bigcap _{
h\in H_v(B)} \Ker (h\otimes _RQ)(t\otimes Q)=\bigcap _{ ht\in
H_{vt}(B)} \Ker (ht\otimes _RQ)=\bigcap _{ h\in H_{vt}(B)} \Ker
(h\otimes _RQ).\]

(2) As $H_{mv}(B)\subseteq H_v(B)$,
\[\bigcap _{
h\in H_v(B)} \Ker (h\otimes _RQ)\subseteq \bigcap _{ h\in
H_{mv}(B)} \Ker (h\otimes _RQ).\] Hence, if $v$ $B$-dominates $u$
with respect to $Q$ we deduce that also $mv$ $B$-dominates $u$
with respect to $Q$.
\end{Proof}

We recall the following property of direct limits.

\begin{Lemma} \label{limitS} Let $M$ be a right $R$-module, and let $\mathcal{S}$ be a
class of finitely presented modules. Then $M\in \varinjlim
\mathcal{S}$ if and only if for any finitely presented module $F$
and any map $u\colon F\to M$ there exists $S\in \mathcal{S}$ and
$v\colon F\to S$ such that $u$ factors through $v$.
\end{Lemma}

\begin{Proof} Assume $M=\varinjlim S_\gamma$ where $(S_{\gamma},u_{\delta\,\gamma})_{\delta\,\gamma\in
I}$. Let $F$ be a finitely presented module and $u\in \Hom
_R(F,M)$. Since $\Hom _R(F,M)$ is canonically isomorphic to
$\varinjlim \Hom _R(F,S_\gamma)$, there exist $\gamma \in I$ and
$v\colon F\to S_\gamma$ such that $u=u_\gamma v$ where
$u_\gamma\colon S_\gamma \to M$ denotes the canonical morphism.

To prove the converse, write $M=\varinjlim F_\alpha$ where
$(F_\alpha,u_{\beta \,\alpha})_{\beta \,\alpha \in I  }$ is a
direct system of finitely presented right $R$-modules. By
hypothesis, for each $\alpha \in I$ there exists $S_\alpha \in
\mathcal{S}$, $v_\alpha \colon F_\alpha \to S_\alpha$ and
$t_\alpha\colon S_\alpha \to M$ such that the canonical map
$u_\alpha \colon F_\alpha \to M$ satisfies $u_\alpha= t_\alpha
v_\alpha$.

Fix $\alpha \in I$. As $\Hom _R(S_\alpha,\varinjlim F_\gamma )$ is
canonically isomorphic to $\varinjlim \Hom _R(S_\alpha,  F_\gamma
)$, there exists $\beta \ge \alpha$ and a commutative diagram:
\[\begin{array}{ccc}
S_{\alpha} & {\stackrel{v'_{\beta\,\alpha}}{\longrightarrow}}&F_\beta\\
t_{\alpha} \downarrow \phantom{t_{\alpha  } }&{\phantom{{u_\beta}}{\swarrow}}{{}_{u_\beta}}
 & \\
M&&
\end{array}\]
Set $u'_{\beta\,\alpha}=v_\beta v'_{\beta\,\alpha}$.

It is not difficult to see that $(S_\alpha,u'_{\beta
\,\alpha})_{\beta \,\alpha \in I  }$ is a direct system of
modules in $\mathcal{S}$ such that $M=\varinjlim S_\alpha$.
\end{Proof}

The next result will provide us with a tool for comparing the relative Mittag-Leffler conditions.
 In fact, we will see in Theorem \ref{2.1.4} that the
   $\B$-stationary
   modules are the modules satisfying the equivalent conditions in \ref{dominatinglimits} for every $B\in\B$ and every  $Q\in R \LMod$,
   while the
   $\Q$-Mittag-Leffler modules are the modules satisfying the equivalent conditions in \ref{dominatinglimits}
   for every $B\in\Mod R $ and   every $Q\in\Q$, see Theorem \ref{AF2}.

\begin{Prop}\label{dominatinglimits} Let $B$ be a right $R$-module, let $Q$ be a  left $R$-module, and let $\mathcal{S}$ be
a class of finitely  presented right $R$-modules. For a right
$R$-module $M\in \varinjlim \mathcal{S}$ the following statements are equivalent.

\noindent{\textrm{(1).}} There
is a direct system of  finitely presented right $R$-modules
$(F_\alpha,u_{\beta \,\alpha})_{\beta \,\alpha \in I  }$    with
$M=\varinjlim (F_\alpha,u_{\beta \alpha})_{\beta , \alpha \in I }$
having the property that for  any $\alpha \in I  $ there exists
$\beta \ge \alpha$ such that $u_{\beta \alpha}$ $B$-dominates the
canonical map $u_\alpha \colon  F_\alpha \to M$ with respect to
$Q$.

\noindent{\textrm{(2).}} Every
direct system of  finitely presented right $R$-modules
$(F_\alpha,u_{\beta \,\alpha})_{\beta \,\alpha \in I  }$    with
$M=\varinjlim (F_\alpha,u_{\beta \alpha})_{\beta , \alpha \in I }$
has the property that for  any $\alpha \in I  $ there exists
$\beta \ge \alpha$ such that $u_{\beta \alpha}$ $B$-dominates the
canonical map $u_\alpha \colon  F_\alpha \to M$ with respect to
$Q$.

\noindent{\textrm{(3).}} For any finitely presented module $F$ (belonging to $\mathcal{S}$) and
any homomorphism $u\colon  F\to M$ there exist a   module $S\in
\mathcal{S}$ and a homomorphism $v\colon  F\to S$ such that
 $u$ factors through $v$,
and $v$ $B$-dominates $u$ with respect to $Q$.

\end{Prop}

\begin{Proof}
$(1)\Rightarrow (3).$ Let $F$ be a finitely presented module and
$u\colon  F\to M$ a homomorphism. Since $\Hom _R(F,M)$ is
canonically isomorphic to $\varinjlim \Hom _R(F,F_\alpha)$, there
exists $\alpha _0\in I  $ and $t\colon  F\to F_{\alpha_0}$ such
that the diagram
\[\begin{array}{ccc}
F & {\stackrel{t}{\longrightarrow}}&F_{\alpha_0}\\
{\scriptstyle u}  \downarrow \phantom{u }&{\phantom{{u_{\alpha_0}}}{\swarrow}}{{}_{u_{\alpha_0}}} & \\
M&&
\end{array}\]
is commutative. By assumption
 there exists
$\beta \ge \alpha _0$ such that $u_{\beta \alpha _0}$
$B$-dominates $u_{\alpha _0}$ with respect to $Q$. Set
$v'=u_{\beta \alpha _0}t$. As $u= u_\beta u_{\beta \alpha _0}t
=u_\beta v'$, we have $u\in H_{v'}(M)$. Moreover, since  $u_{\beta
\alpha _0}$ $B$-dominates $u_{\alpha _0}$, it follows from
Lemma~\ref{composition}(1) that $v'=u_{\beta \,\alpha _0}t$
$B$-dominates $u=u_{\alpha _0}t$ with respect to $Q$.

By hypothesis, $M=\varinjlim S_\gamma$ for a directed system
$(S_\gamma,u'_{\delta \,\gamma})_{\delta \,\gamma \in J}$ of
modules in $\mathcal{S}$. As $F_{\beta}$ is finitely presented,
there exist $\gamma$ in $J$ and $m\colon F_\beta \to S_\gamma$
such that the diagram
\[\begin{array}{ccc}
F_{\beta} & {\stackrel{m}{\longrightarrow}}&S_{\gamma}\\
{\scriptstyle u_{\beta} }\downarrow \phantom{u_{\beta} }&{\phantom{{u\,'_\gamma}}{\swarrow}}{{}_{u\,'_\gamma}} & \\
M&&
\end{array}\]
commutes. Set $v=m\,v'$. Then $u$ factors through $v$ and, by
Lemma~\ref{composition}(2), $v$ $B$-dominates $u$ with respect to
$Q$.

$(3)\Rightarrow (2).$ Consider a  direct system of  finitely presented right $R$-modules
$(F_\alpha,u_{\beta \,\alpha})_{\beta \,\alpha \in I  }$    with $M=\varinjlim (F_\alpha,u_{\beta
\alpha})_{\beta , \alpha \in I  }$. Fix $\alpha _0 \in I  $. We have to verify the existence of
$\beta\ge \alpha_0$ such that $u_{\beta \alpha_0}$ $B$-dominates $u_{\alpha_0}$ with respect to
$Q$.   Applying the hypothesis with $u=u_{\alpha _0} \colon F_{\alpha_0} \to M$ we deduce that
there exist a   module $S\in \mathcal{S}$, $v\colon F_{\alpha_0}\to S $ and $t\colon  S\to M$ such
that the diagram
\[\begin{array}{ccc}
F_{\alpha_0} & {\stackrel{v}{\longrightarrow}}&S\\
{\scriptstyle u_{\alpha_0} } \downarrow \phantom{u_{\alpha _0} }&{\phantom{{t}}{\swarrow}}{{}_{t}} & \\
M&&
\end{array}\]
is commutative and $v$ $B$-dominates $u_{\alpha_0}$ with respect
to $Q$. As $S$ is finitely presented, $\Hom _R(S,M)$ is canonically
isomorphic to $\varinjlim \Hom _R(S,F_\alpha)$. Hence there exist
$\beta '\ge \alpha _0 $ and $t'\colon  S\to F_{\beta'}$ such that
the diagram
\[\begin{array}{ccc}
F_{\alpha_0} & {\stackrel{v}{\longrightarrow}}&S\\
{\scriptstyle u_{\alpha_0} }\downarrow \phantom{u_{\alpha _0} }&{\phantom{{t}}{\swarrow}}{{}_{t}} & \phantom{t'}\downarrow t'\\
M&{\stackrel{u_{\beta'}}{\longleftarrow}}& F_{\beta '}
\end{array}\]
is commutative. Since  $u_{\beta '}u_{\beta '\, \alpha_0}=u_{\beta
'}t'v$, there exists $\beta \ge \beta'$ such that $u_{\beta
\beta'}u_{\beta '\alpha_0}=u_{\beta \beta'}t'v$, that is,
$u_{\beta \alpha_0}=m\,v$ where $m=u_{\beta \beta'}t'$. By
Lemma~\ref{composition}(2), $u_{\beta \alpha_0}$ $B$-dominates
$u_{\alpha_0}$ with respect to $Q$.

Similarly, to see that condition $(3)$ restricted to modules $F$ belonging to $\mathcal S$ implies  $(1)$, we proceed as in
$(3)\Rightarrow (2)$ but considering a direct system of  finitely presented right $R$-modules
$(F_\alpha,u_{\beta \,\alpha})_{\beta \,\alpha \in I  }$    with $M=\varinjlim (F_\alpha,u_{\beta
\alpha})_{\beta , \alpha \in I  }$ such that all $F_\alpha \in \mathcal{S}$.
\end{Proof}

Observe that the condition   $M\in \varinjlim \mathcal{S}$ in the hypothesis of  Proposition \ref{dominatinglimits} is also necessary.
 This can be deduced from condition (3)
by employing Lemma~\ref{limitS}.

 \medskip

We will need the following result.

\begin{Prop}{\rm \cite[Proposition~2.1.1]{RG}} \label{trick}
Let $u\colon  M\to N$ and $h\colon  M\to B$ be right $R$-module
homomorphisms. The following statements are equivalent.
\begin{itemize}
\item[(i)] $\Ker (u\otimes _RQ)\subseteq \Ker (h\otimes Q)$ for all left $R$-modules $Q$.
\item[(ii)] $\Ker (u\otimes _R B^\ast)\subseteq \Ker (h\otimes B^\ast)$.
\end{itemize}
If $\mathrm{coker}(u)$ is finitely presented, the following statement is further equivalent.
\begin{itemize}
\item[(iii)] $h$ factors through $u$.
\end{itemize}
\end{Prop}

We can now interpret the property ``$B$-dominates'' in terms of $H$-subgroups.

\begin{Prop} \label{2.1.1} Let $B$ be a  right
$R$-module. Let $u\colon  M\to N$ and $v\colon  M\to M'$ be right
$R$-module homomorphisms.
 If $H_v(B)\subseteq H_u(B)$, then
$v$ $B$-dominates $u$. The converse implication holds true provided $\mathrm{coker}(u)$ is finitely presented.
\end{Prop}

\begin{Proof}  Let  $h\colon  M\to B\in H_v(B)$.
By hypothesis, there exists $h'\colon  N\to B$ such that $h=h'u$.
Hence, for any left $R$-module $Q$
\[\Ker (u\otimes _RQ)\subseteq \Ker (h'\otimes Q)(u\otimes _RQ)\subseteq\Ker (h\otimes _RQ). \]
This shows the claim.

For the converse implication, assume that $v$ $B$-dominates $u$ and $\mathrm{coker}(u)$ is finitely presented.
Let  $h\in
H_v(B)$. Then $\Ker (u\otimes
_RQ)\subseteq\Ker (h\otimes _RQ)$ for any left $R$-module $Q$.  By
\ref{trick} this means that  $h\in H_u(B)$.
\end{Proof}

\begin{Lemma}\label{directlimits} Let   $B$ be a  right $R$-module and let
$Q$ be a left $R$-module. Let further $(F_\alpha,u_{\beta \,\alpha})_{\beta \,\alpha \in
I  }$ be a direct system of  finitely presented right $R$-modules
with $M=\varinjlim (F_\alpha,u_{\beta \alpha})_{\beta , \alpha \in
I  }$.   For $\alpha$, $\beta \in I  $
with $\beta\ge\alpha$,  the following statements hold true.
\begin{itemize}
\item[(1)] $u_{\beta \alpha}$ $B$-dominates the   canonical map $u_\alpha \colon
F_\alpha \to M$  with respect to $Q$ if and only if $u_{\beta \alpha}$ $B$-dominates $u_{\gamma
\alpha}$ with respect to $Q$ for any $\gamma \ge \alpha$.
\item[(2)] $u_{\beta \alpha}$ $B$-dominates the   canonical map $u_\alpha \colon
F_\alpha \to M$  if and only if  $H_{u_{\beta \alpha}}(B)=\bigcap _{\gamma \ge \alpha}H_{u_{\gamma
\alpha}}(B)$.
\end{itemize}
\end{Lemma}
\begin{Proof}
(1) To show  the only-if-part, fix $\gamma \ge \alpha$. As $u_\alpha =u_\gamma u_{\gamma
\alpha}$, for any left $R$-module $Q$ \[\Ker (u_{\gamma
\alpha}\otimes _RQ)\subseteq\Ker (u_\alpha\otimes _RQ)\subseteq
\bigcap _{ h\in H_{u_{\beta \alpha}}(B ) } \Ker (h\otimes _RQ).\]
Therefore $u_{\beta \alpha}$ $B$-dominates $u_{\gamma \alpha}$
with respect to $Q$.
\\
The converse implication is clear from the properties of direct
limits.
\\
(2) By (1), $u_{\beta \alpha}$ $B$-dominates $u_\alpha$ if and only if $u_{\beta \alpha}$ $B$-dominates $u_{\gamma\alpha}$ for any $\gamma \ge \alpha$.
As $\mathrm{coker}\, u_{\gamma \alpha}$ is finitely presented, we know from
Proposition~\ref{2.1.1}   that the latter is equivalent to
 $H_{u_{\beta \alpha}}(B)\subseteq H_{u_{\gamma
\alpha}}(B)$ for any $\gamma \ge \alpha$.
But this   means $H_{u_{\beta \alpha}}(B)=\bigcap _{\gamma \ge \alpha}H_{u_{\gamma
\alpha}}(B)$ by Lemma \ref{alphabeta}.
\end{Proof}

{}From Lemma \ref{directlimits} and Definition \ref{Bstat}, we immediately obtain the announced
characterization of $B$-stationary modules.

\begin{Th}\label{2.1.4}
Let  $B$ be a right $R$-module, and let $\mathcal{S}$ be a
class of finitely  presented modules. For a right $R$-module $M\in \varinjlim \mathcal{S}$, the following statements are
equivalent.

\noindent{\textrm{(1).}}   $M$ is
$B$-stationary.

\noindent{\textrm{(2).}} There
is a direct system of  finitely presented right $R$-modules
$(F_\alpha,u_{\beta \,\alpha})_{\beta \,\alpha \in I  }$    with
$M=\varinjlim (F_\alpha,u_{\beta \alpha})_{\beta , \alpha \in I }$
having the property that for  any $\alpha \in I  $ there exists
$\beta \ge \alpha$ such that $u_{\beta \alpha}$ $B$-dominates the
canonical map $u_\alpha \colon  F_\alpha \to M$.

\noindent{\textrm{(3).}} For any finitely presented module $F$ (belonging to $\mathcal{S}$) and
any homomorphism $u\colon  F\to M$ there exist a   module $S\in
\mathcal{S}$ and a homomorphism $v\colon  F\to S$ such that
 $u$ factors through $v$,
and $v$ $B$-dominates $u$.

\end{Th}

  We close this section with some closure properties of the class $\B$ in the definition of ``$\B$-dominating''. Let us first prove the following preliminary result.
\begin{Prop}\label{fpres} Let $B$ be a right $R$-module. Let $u\colon  M\to N$ and $v\colon  M\to M'$ be right
$R$-module homomorphisms. Then the following statements are
equivalent.
\begin{itemize}
\item[(1)] $v$ $B$-dominates $u$.
\item[(2)] For  any finitely presented left $R$-module $Q$
\[\Ker (u\otimes _RQ)\subseteq \bigcap _{
h\in H_v(B)} \Ker (h\otimes _RQ). \]
\item[(3)] For  any (finitely presented) left $R$-module $Q$
\[\Ker (u\otimes _RQ)\subseteq  \Ker (\widetilde{h}\otimes _RQ). \]
where $\widetilde{h}\colon M\to \prod _{H_v(B)} B$ is the product
map induced by all $h\in H_v(B)$.
\end{itemize}
\end{Prop}

\begin{Proof} We follow the idea in the proof of
\cite[Proposition~2.1.1]{RG}.

Fix  $h\in H_v(B)$. Consider the push-out diagram
\[\begin{array}{ccc}
M&\stackrel{h}\to&B\\
{\scriptstyle u} \downarrow \phantom{u}&&\phantom{u'}\downarrow {u'}\\
N&\stackrel{h'}\to&N'
\end{array}\]
Recall that it will stay a push-out diagram when we apply the
functor $-\otimes _RQ$ for any left module $Q$. Hence we have the
exact sequence
\[0\to \Ker(u\otimes _RQ)\cap \Ker(h\otimes _RQ)\to \Ker(u\otimes _RQ) \stackrel{h\otimes _RQ}\to \Ker(u'\otimes _RQ)\to 0.\]
This shows that, for any left module $Q$, $\Ker(u\otimes
_RQ)\subseteq \Ker(h\otimes _RQ)$ if and only if $\Ker(u'\otimes
_RQ)=0$, that is, if and only if $u'$ is a pure monomorphism.

Since a morphism is a pure monomorphism if and only if it is  a monomorphism when tensoring by finitely
presented modules, we deduce that $(1)$ and $(2)$ are equivalent
statements.

To prove that $(2)$ and $(3)$ are equivalent, note that
$\bigcap _{ h\in H_v(B)} \Ker (h\otimes _RQ)$ is the kernel of the
product map induced by all homomorphisms $h\otimes _RQ$ with $h\in
H_v(B)$. When  $Q$ is finitely presented, the natural morphism
$\rho\colon  \prod _{H_v(B)}B\otimes _RQ\to \prod
_{H_v(B)}(B\otimes _RQ)$ is an isomorphism. Hence
\[\bigcap _{
h\in H_v(B)} \Ker (h\otimes _RQ)=\Ker(\widetilde{h}\otimes
_RQ)\]
and the statement is verified.
To obtain the statement for arbitrary $Q$, proceed
as in the proof of $(1)\Leftrightarrow (2)$.
\end{Proof}

\begin{Prop} Let $u\colon  M\to N$ and $v\colon  M\to M'$ be right
$R$-module homomorphisms. Let $\B$ be a class of right $R$-modules
such that  $v$ $\B$-dominates $u$. Then
\begin{itemize}
\item[(i)]
$v$ $\B'$-dominates $u$, where $\B'$ denotes the class of all pure submodules of modules in $\B$.
\item[(ii)] $v$ {\rm Prod} $\B$-dominates
$u$.
\item[(iii)]  $v$ $\varinjlim\B$-dominates $u$ provided that   $M'$ is finitely presented.
\end{itemize}
\end{Prop}

\begin{Proof} (i) Let $B\in\B$, and assume that  the inclusion $\epsilon \colon  C\to B$ is a
pure monomorphism. If $h\in H_v(C)$, then $\epsilon h\in H_v(B)$, and  $\Ker
(h\otimes Q)= \Ker (\epsilon h\otimes Q)$ contains $\Ker (u\otimes Q)$.

(ii)  By (i) it is enough to consider modules of the form
$\prod_{i\in I}B_i$ where $\{B_i\}_{i\in I}$ is a family of
modules in $\B$. Let $Q$ be a finitely presented module. As the
canonical morphism $\rho\colon  (\prod _{i\in I}B_i)\otimes Q\to
\prod _{i\in I}B_i\otimes _RQ$ is an isomorphism and, as any $h\in
H_v(\prod _{i\in I}B_i)$ is induced by a family $(h_i)_{i\in I}$
where $h_i \in H_v(B_i)$ for any $i\in I$, we deduce that
\[\bigcap _{
h\in H_v(\prod _{i\in I}B_i)} \Ker (h\otimes _RQ)=\bigcap
_{h\in H_v(B_i), \, i\in I} \Ker
(h\otimes _RQ).\] Then the claim follows from \ref{fpres}.

(iii) Let  $\{B_i, f_{ji}\}_{i\in I}$ be a direct system of
modules in $\B$, and let $h\in H_v(\varinjlim B_i)$. Then
$h=h'\,v$ for some $f\colon M'\to \varinjlim B_i$.
 As $M'$ is finitely presented, there exists $j\in I$ such that $h'$ factors
through the canonical map $f_j\colon  B_j\to \varinjlim B_i$. So,
there exists $g\colon  M'\to B_j$ such that $h'=f_jg$, thus
$h=f_j\,g\,v$ with $g\,v\in H_v(B_j)$. Hence, for any left module
$Q$, we have $\Ker (u\otimes _RQ)\subseteq \Ker (gv\otimes
_RQ)\subseteq \Ker (h\otimes _RQ)$.
\end{Proof}


\section{$\Q$-Mittag-Leffler modules revisited}\label{revisit}

As a next step towards establishing a relationship between
 $\Q$-Mittag-Leffler and $\B$-stationary modules, we provide
 a characterization of $\Q$-Mittag-Leffler modules in terms of dominating maps. It is inspired to work of Azumaya and Facchini \cite[Theorem 6]{AF}.

\begin{Th} \label{AF2} Let $\Q$ be a class of left $R$-modules, and let $\mathcal{S}$ be
a class of finitely  presented right $R$-modules. For a right $R$-module $M\in \varinjlim \mathcal{S}$, consider the
following statements.
\begin{itemize}
\item[(1)] $M$   is  $\Q$-Mittag-Leffler.
\item[(2)] Every
direct system of  finitely presented right $R$-modules $(F_\alpha,u_{\beta \,\alpha})_{\beta
\,\alpha \in I  }$    with $M=\varinjlim (F_\alpha,u_{\beta \alpha})_{\beta , \alpha \in I }$ has
the property that for  any $\alpha \in I  $ there exists $\beta \ge \alpha$ such that $u_{\beta
\alpha}$ dominates the canonical map $u_\alpha \colon  F_\alpha \to M$ with respect to $\Q$.

\item[(3)] For every finitely presented module $F$ (belonging to $\mathcal S$) and every homomorphism $u\colon F\to M$ there
are a   module $S\in \mathcal{S}$ and a homomorphism $v\colon F\to S$ such that
$u$ factors through $v$ and  $\Ker (u \otimes_R Q ) = \Ker (v\otimes _RQ) $ for all $Q\in\Q$.

\item[(4)] For every countable (finite) subset $X$ of $M$ there are a countably presented
$\Q$-Mittag-Leffler module $N\in \varinjlim \Scal$ and a homomorphism $v\colon N\to M$ such that
$X\subseteq v(N)$ and $ v\otimes _RQ$ is a monomorphism for all $Q\in\Q$.

\item[(5)]  For every finitely generated submodule $M_0$ of $ M$ there are a
finitely presented module $S\in \mathcal{S}$ and a homomorphism $w\colon M_0\to S$ such that the
embedding $\epsilon\colon M_0\to M$ factors through $w$ and  $\Ker (\epsilon \otimes_R Q ) = \Ker
(w\otimes _RQ) $ for all $Q\in\Q$.

\end{itemize}
Then $(1)$, $(2)$, $(3)$ and $(4)$ are equivalent statements, and $(5)$  implies the other
statements. Moreover, if $R\in \Q$, then all statements are equivalent.
\end{Th}

\begin{Proof}
$(1)\Rightarrow (2)$. Assume for a contradiction that there is a direct system of
finitely presented modules $(F_\alpha,u_{\beta \,\alpha})_{\beta \,\alpha \in I  }$    with
$M=\varinjlim (F_\alpha,u_{\beta \alpha})_{\beta , \alpha \in I  }$ such that there exists $\beta
\in I$ satisfying that for any $\beta \ge \alpha $ there exists $Q_{\beta}\in \Q$ such that
\[\mathrm{ker}(u_\alpha \otimes Q_{\beta})\nsubseteqq \mathrm{ker}(u_{\beta\,\alpha} \otimes Q_{\beta}).\]
If $x_1,\dots, x_n$ is a generating set of $F_\alpha$, then for each $\beta \ge \alpha $ we can choose
\[a_\beta=\sum _{i=1}^nx_i\otimes q_\beta^i\in \mathrm{ker}(u_\alpha \otimes Q_{\beta})\setminus
\mathrm{ker}(u_{\beta\,\alpha} \otimes Q_{\beta}).\] Set $x=\sum _{i=1}^nx_i\otimes
(q_\beta^i)_{\beta \ge \alpha} \in F_\alpha \otimes \prod_{\beta \ge \alpha} Q_\beta$. Consider
the commutative diagram,
$$\begin{array}{ccccc}&&F_\alpha \otimes \prod_{\beta \ge \alpha} Q_\beta & {\stackrel{u_\alpha \otimes \prod_{\beta \ge \alpha} Q_\beta
}{\longrightarrow}} & M\otimes \prod_{\beta \ge \alpha} Q_\beta\\
&&\rho'\downarrow \phantom{\rho'}& & \rho  \downarrow \phantom{\rho }\\
&& \prod_{\beta \ge \alpha} (F_\alpha \otimes  Q_\beta)&{\stackrel{\prod_{\beta \ge \alpha}
(u_\alpha \otimes  Q_\beta)}{\longrightarrow}} & \prod_{\beta \ge \alpha} (M \otimes
Q_\beta)\end{array}
$$
As $\left( \prod_{\beta \ge \alpha} (u_\alpha \otimes  Q_\beta) \right)\rho'(x)=0$ and, by
hypothesis, $\rho$ is injective, we deduce that $(u_\alpha \otimes \prod_{\beta \ge \alpha}
Q_\beta)(x)=0$. Since $M\otimes \prod_{\beta \ge \alpha} Q_\beta=\varinjlim \left( F_\gamma
\otimes \prod_{\beta \ge \alpha} Q_\beta \right)$, there exists $\beta _0\ge \alpha$ such that
$x\in \mathrm{ker} (u_{\beta _0\,\alpha}\otimes \prod_{\beta \ge \alpha} Q_\beta)$. The
commutativity of the diagram
$$\begin{array}{ccccc}&&F_\alpha \otimes \prod_{\beta \ge \alpha} Q_\beta & {\stackrel{u_{\beta _0\,\alpha} \otimes \prod_{\beta \ge \alpha} Q_\beta
}{\longrightarrow}} & F_{\beta _0}\otimes \prod_{\beta \ge \alpha} Q_\beta\\
&&\cong\downarrow \phantom{\cong}& & \cong  \downarrow \phantom{\cong }\\
&& \prod_{\beta \ge \alpha} (F_\alpha \otimes  Q_\beta)&{\stackrel{\prod_{\beta \ge \alpha}
(u_{\beta _0\,\alpha} \otimes  Q_\beta)}{\longrightarrow}} & \prod_{\beta \ge \alpha} (F_{\beta
_0} \otimes Q_\beta)\end{array}
$$
implies that, for any $\beta \ge \alpha$, $a_\beta \in \mathrm{ker}(u_{\beta _0\,\alpha} \otimes
Q_\beta)$. In particular, $a_{\beta _0}\in \mathrm{ker}(u_{\beta _0\,\alpha} \otimes Q_\beta)$
which is a contradiction.

$(2)\Rightarrow (1)$.  Fix a direct system of finitely presented right $R$-modules
$(F_\alpha,u_{\beta \,\alpha})_{\beta \,\alpha \in I  }$ with $M=\varinjlim (F_\alpha,u_{\beta
\alpha})_{\beta , \alpha \in I }$.

Let $\{Q_k\}_{k\in K}$ be a family of modules of $\Q$, and let $x\in \mathrm{ker} \rho$ where
$\rho \colon M\otimes \prod _{k\in K}Q_k\to \prod _{k\in K} (M\otimes Q_k)$ denotes the natural
map. Since $M\otimes \prod _{k\in K}Q_k=\varinjlim \left( F_\alpha \otimes \prod _{k\in K}Q_k
\right)$ there exists $\alpha \in I$ and $x_\alpha =\sum _{i=1}^nx_i\otimes (q_k^i)_{k\in K}\in
F_\alpha \otimes \prod _{k\in K}Q_k$ such that $x=(u_\alpha \otimes \prod _{k\in
K}Q_k)(x_\alpha)$. The commutativity of the diagram
$$\begin{array}{ccccc}&&F_\alpha \otimes \prod _{k\in K}Q_k & {\stackrel{u_\alpha \otimes \prod _{k\in K}Q_k}
{\longrightarrow}} & M\otimes \prod _{k\in K}Q_k\\
&&\rho'\downarrow \phantom{\rho'}& & \rho  \downarrow \phantom{\rho }\\
&& \prod_{k\in K} (F_\alpha \otimes  Q_k)&{\stackrel{\prod_{k\in K} (u_\alpha \otimes
Q_k)}{\longrightarrow}} & \prod_{k\in K} (M \otimes Q_k)\end{array}
$$
implies that, for each $k\in K$, $\sum _{i=1}^nx_i\otimes q_k^i\in
\mathrm{ker}(u_\alpha \otimes Q_k)$.

Let $\beta \ge \alpha$ be such that $u_{\beta \alpha}$ dominates the canonical map $u_\alpha $
with respect to $\Q$. The commutativity of the diagram
$$\begin{array}{ccccc}&&F_\alpha \otimes \prod_{k\in K} Q_k & {\stackrel{u_{\beta \,\alpha} \otimes \prod_{k\in K}
Q_k
}{\longrightarrow}} & F_{\beta }\otimes \prod_{k\in K} Q_k\\
&&\cong\downarrow \phantom{\cong}& & \cong  \downarrow \phantom{\cong }\\
&& \prod_{k\in K} (F_\alpha \otimes  Q_k)&{\stackrel{\prod_{k\in K} (u_{\beta \,\alpha} \otimes
Q_k)}{\longrightarrow}} & \prod_{k\in K} (F_{\beta } \otimes Q_k)\end{array}
$$
implies that $(u_{\beta \,\alpha} \otimes \prod_{k\in K} Q_k)(x_\alpha)=0$. Hence $x=(u_\beta
u_{\beta \,\alpha} \otimes \prod_{k\in K} Q_k) (x_\alpha)=0$.

By Proposition~\ref{dominatinglimits}, we already know that $(2)$ and $(3)$ are equivalent
statements.

$(3)\Rightarrow (4)$.
Let $X=\{x_1,x_2,\dots \}\subseteq M$. We shall construct inductively a countable direct system
$(S_n,f_n\colon S_n\to S_{n+1})_{n\ge 0}$ of modules in $\mathcal{S}$  and a sequence of maps
$(v_n\colon S_n\to M)_{n\ge 0}$  such that $v_n=v_{n+1}f_n$ and $\{x_1,\dots ,x_n\}\subseteq
v_n(S_n)$.

Set $S_0=0$ and let $v_0$ be the zero map. Let $n\ge 0$ and assume as inductive
hypothesis that $S_m$ and $v_m$ have been constructed for any $m\le n$.
Let $u\colon S_n\oplus R\to M$ be defined as $u(g,r)=v_n(g)+x_{n+1}r$ for any $(g,r)\in S_n\oplus
R$. By $(3)$, there exist $S_{n+1}\in \mathcal{S}$, $v\colon S_n\oplus R\to S_{n+1}$, and
$v_{n+1}\colon S_{n+1}\to M$ such that $u=v_{n+1}v$ and $\mathrm{ker} (v\otimes Q)=\mathrm{ker}
(u\otimes Q)$ for all $Q\in \Q$. Let $\varepsilon \colon S_n\to S_n\oplus R$ denote the canonical
inclusion and set $f_{n}=v\circ \varepsilon.$ Then
$v_n=u\varepsilon=v_{n+1}(v\,\varepsilon)=v_{n+1}f_n$.
This completes the induction
step.
Note moreover that also $\mathrm{ker}
(v_n\otimes Q)=\mathrm{ker} (f_{n}\otimes Q)$ for all $Q\in \Q$.

Set $N=\varinjlim S_n$ and $v=\varinjlim v_n$. Then $N$ is countably presented.
As for any $Q\in \Q$, $\mathrm{ker} (v\otimes
Q)=\varinjlim \mathrm{ker}(v_n\otimes Q)$ and $\mathrm{ker}(v_n\otimes Q)=\mathrm{ker}(f_n\otimes
Q)$, we deduce that $v\otimes Q$ is injective.

To show that $N$ is $\Q$-Mittag-Leffler we verify that $N$ satisfies $(2)$, as we already know that
$(1)$ and $(2)$ are equivalent. By Proposition~\ref{dominatinglimits} it is enough to check the
condition for the direct  system $(S_n,f_n)_{n\ge 0}$ and the canonical maps
$u_n\colon S_n\to M$.

Notice that $vu_n=v_{n}$. Therefore, for any $Q\in \Q$
\[\mathrm{ker}(u_n\otimes Q)\subseteq \mathrm{ker}(v_n\otimes Q)=
\mathrm{ker}(f_n\otimes Q)\] from which we conclude that
$f_n$ dominates $u_n$ with respect to $\Q$.

$(5)\Rightarrow (4)$ is proven similarly.

$(4)\Rightarrow (1).$   Consider a family $(Q_k)_{k\in K}$ in $\Q$,  and  an element $x$ in the
kernel of \linebreak \mbox{$\rho\colon M\otimes \prod _{k\in K}Q_k\to \prod _{k\in
K}(M\bigotimes_RQ_k).$} Then there are a
$\Q$-Mittag-Leffler module $N$ and a homomorphism
$v\colon N\to M$  such that $x$ lies in the image of  $(v\otimes \prod _{k\in K}Q_k)$, and  $
v\otimes _RQ$ is a monomorphism for all $Q\in\Q$. In the commutative diagram
$$\begin{array}{ccccc}&&N \otimes \prod _{k\in K}Q_k & {\stackrel{v\otimes
\prod _{k\in K}Q _k
}{\longrightarrow}} & M\otimes \prod _{k\in K}Q_k\\
&&\rho'\downarrow \phantom{\rho'}& & \rho  \downarrow \phantom{\rho }\\
&& \prod _{k\in K}(N\bigotimes_RQ_k)&{\stackrel{\prod _{k\in K}(v\otimes Q_k)}{\longrightarrow}} &
\prod _{k\in K}(M\bigotimes_RQ_k)\end{array}
$$
we then have that $\rho'$ is injective because $N$ is $\Q$-Mittag-Leffler, and $\prod _{k\in
K}(v\otimes Q_k)$ is injective by assumption on $v$. This shows that $x=0$.

Assume now that $R\in \Q$. To show  $(3)\Rightarrow (5)$, we proceed as in the proof of $(2)\Rightarrow
(3)$ in \cite[Theorem 6]{AF}. We take an epimorphism $p\colon F\to M_0$ from a finitely generated
free module $F$, set $u=\epsilon\, p$, and construct $v$ as in condition (3).  Note that $\Ker u=
\Ker v$ since $\Q$ contains $R$. We  thus obtain $w\colon M_0\to S$ and $t\colon S\to M$ such that
$v=w\,p$ and $\epsilon=t\,w$. To show $\Ker (\epsilon \otimes_R Q ) = \Ker
(w\otimes _RQ) $ for all $Q\in\Q$ it is enough to verify the inclusion $\subseteq$. So, take a left
$R$-module $Q\in\Q$ and $y\in \Ker(\epsilon\otimes Q)$. Note that $y=(p\otimes Q)(x)$ for some $x\in
F\otimes Q$. Then $(u\otimes Q)(x)=(\epsilon\,p\otimes Q)(x)=0$, hence $(w\otimes Q)(y)= (v\otimes Q)(x)=0$.
\end{Proof}

Condition $(4)$ in Theorem~\ref{AF2} gives the following characterization of $\Q$-Mittag-Leffler
modules.

\begin{Cor} \label{QML}Let $\Q$ be a class of left $R$-modules, and let $\mathcal{S}$ be
a class of finitely  presented right $R$-modules. For a fixed right $R$-module $M\in \varinjlim
\mathcal{S}$  denote by $\mathcal{C}$ the class of its countably generated submodules $N$ such
that $N$ is $\Q$-Mittag-Leffler and   the inclusion $N\subseteq M$ remains injective when
tensoring with any module $Q \in \Q$.

Then $M$ is $\Q$-Mittag-Leffler if and only if $M$ is a directed union of modules in
$\mathcal{C}$.

Moreover, if $R\in \Q$, the modules in $\mathcal C$ can be taken  countably presented and in
$\varinjlim \mathcal{S}$.
\end{Cor}

\begin{Proof} For the only-if implication, we follow the notation of Theorem~\ref{AF2} (4). We only have to
prove that $v(N)$ is a $\Q$-Mittag-Leffler  module and that
the inclusion $\varepsilon \colon v(N)\to M$  remains injective when tensoring with any module $Q \in \Q$.
Let $\{Q_k\}_{k\in K}$ be a family of modules in $\Q$. Consider the commutative
diagram,
$$\begin{array}{ccccccc}N \otimes \prod _{k\in K}Q_k & {\stackrel{v\otimes
\prod _{k\in K}Q _k }{\longrightarrow}} & v(N) \otimes \prod _{k\in K}Q_k &
{\stackrel{\varepsilon\otimes \prod _{k\in K}Q _k
}{\longrightarrow}} &M\otimes \prod _{k\in K}Q_k\\
\rho _1\downarrow \phantom{\rho _1}& &\rho _2\downarrow \phantom{\rho _2}& & \rho  \downarrow \phantom{\rho }\\
 \prod _{k\in K}(N\bigotimes_RQ_k)&{\stackrel{\prod _{k\in K}(v\otimes Q_k)}{\longrightarrow}} &
\prod _{k\in K}(v(N)\bigotimes_RQ_k)&{\stackrel{\prod _{k\in K}(\varepsilon\otimes
Q_k)}{\longrightarrow}} & \prod _{k\in K}(M\bigotimes_RQ_k)\end{array}
$$
Note that $v\otimes \prod _{k\in K}Q _k: N \otimes \prod _{k\in K}Q_k\to v(N) \otimes \prod _{k\in K}Q_k$ is surjective,
therefore if $x\in v(N) \otimes \prod
_{k\in K}Q_k$ satisfies $(\varepsilon\otimes \prod _{k\in K}Q _k)(x)=0$, then there exists $y\in
N \otimes \prod _{k\in K}Q_k$ such that $x=(v\otimes \prod _{k\in K}Q _k) (y)$ and $(\varepsilon v\otimes \prod _{k\in K}Q _k)(y)=0$. Since
  $\rho\, (\varepsilon v \otimes \prod _{k\in K}Q _k) = (\prod _{k\in K}(v\otimes Q_k))\, \rho _1$ is an injective map,
 we infer $y=0$, so $x=0$. This shows that $\varepsilon\otimes \prod
_{k\in K}Q _k$ is injective. Then also $\left(\prod _{k\in K}(\varepsilon\otimes Q_k)\right)\,\rho _2= \rho\,
(\varepsilon\otimes \prod _{k\in K}Q _k) $ is injective, and so is $\rho _2$.

To prove the converse implication proceed as in the proof of $(4)\Rightarrow (1)$ of
Theorem~\ref{AF2}.

The statement for the case when $R\in \Q$ is clear because then the map $v$ in Theorem~\ref{AF2} is injective,
so $N$ is isomorphic to $v(N)$.
\end{Proof}

\begin{Cor}\label{MLcountable} Let $\Q$ be a class of left $R$-modules containing $R$.
Then every countably generated $\Q$-Mittag-Leffler
 right $R$-module is countably presented.
\end{Cor}

\medskip

Now we can start  relating $\Q$-Mittag-Leffler and $\B$-stationary modules.

\begin{Lemma} \label{zero}
Let  $\Q$ be a class of left $R$-modules, and let $B, M$ be right $R$-modules. Assume that $M$ is  $B$-stationary.
Write $M=\varinjlim F_\alpha$ where $(F_\alpha,u_{\beta \,\alpha})_{\beta \,\alpha \in
I  }$ is a direct system of finitely presented modules.
If for all $\alpha,\beta\in I$ with  $\beta\ge \alpha$ and  all $Q\in\Q$
$$\Ker (u_{\beta\alpha} \otimes_R Q ) =\bigcap _{
h\in H_{u_{\beta\alpha}}(B)} \Ker (h\otimes _RQ) $$
then $M$ is a $\Q$-Mittag-Leffler module.
\end{Lemma}
\begin{Proof}
 Fix $\alpha \in I  $, and denote
by $u_\alpha \colon  F_\alpha \to M$ the canonical map. As $M$ is $B$-stationary,  we infer from
Proposition \ref{2.1.4} that
 there exists
$\beta \ge \alpha $ such that $u_{\beta \alpha}$
$B$-dominates the canonical map $u_{\alpha}$, that is
$$\Ker (u_\alpha\otimes _RQ)\subseteq \bigcap _{
h\in H_{u_{\beta\alpha}}(B)} \Ker (h\otimes _RQ) $$ for all  left $R$-modules $Q$. Our assumption
implies that $\Ker (u_\alpha\otimes _RQ)\subseteq \Ker (u_{\beta\alpha} \otimes_R Q )$ for all $Q\in\Q$, so
Theorem~\ref{AF2} gives the desired conclusion.
\end{Proof}

Before we continue our discussion of the general case, let us
 notice the following projectivity criteria for countably generated flat modules that improves
\cite[Propostion~2.5]{abh}, and clarifies the proof of \cite[Theorem~2.2]{drinfeld}.

\begin{Cor}\label{cp} Let $M$ be a countably generated right flat module. Then $M$ is projective if and only
if $M$ is $R$-stationary.
\end{Cor}

\begin{Proof} To see that a countably generated projective module is $R$-stationary use for example
Theorem~\ref{BH}.

Assume  that $M$ is countably generated, flat and $R$-stationary. Then $M$ is also $R^{(\N)}$-stationary by \ref{closureB}. Let
$(F_{\alpha}, u_{\beta \alpha})_{\alpha \, \beta \in I}$ be a direct system of finitely generated
free modules such that $M=\varinjlim F_\alpha$. Notice that for each $\beta \in I$ we have a
split monomorphism $t_\beta \colon F_\beta \to R^{(\N)}$, hence $t_\beta\otimes Q$ is a split monomorphism for any
left $R$-module $Q$. This implies that the criterion of Lemma~\ref{zero} is fulfilled for any left
$R$-module, hence $M$ is a Mittag-Leffler module. Now we can conclude either by using
\cite[
2.2.2]{RG} or arguing that then $M$ is $R$-Mittag-Leffler,
hence countably presented by Corollary \ref{cp}, and then use \cite[Propostion~2.5]{abh}.
\end{Proof}

\begin{Ex}\label{injective-flat}{\rm \cite{goodearl}
If $\mathcal{Q}$ denotes the class of flat left
$R$-modules,  condition
$(5)$ in Theorem~\ref{AF2} is equivalent to:

\begin{quote} (5')  For any finitely generated submodule $M_0$ of $M$
there are a finitely presented module $S$ and a homomorphism
$w\colon M_0\to S$ such that the embedding $\epsilon\colon M_0\to
M$ factors through $w$.
\end{quote}

\noindent
We thus recover  a characterization due to Goodearl of the modules that are Mittag-Leffler
with respect to the class of flat modules
\cite[Theorem~1]{goodearl}.
In particular, if $R$ is  right noetherian, then $(5')$ is trivially satisfied, and so any right $R$-module is
Mittag-Leffler with respect to the class of flat modules (cf.
\cite{goodearl}).
}
\end{Ex}

\section{Relating $\B$-stationary and $\Q$-Mittag-Leffler modules}\label{relat}

Throughout this section, we fix  a right $R$-module $M$ together with    a
direct system of finitely presented modules
$(F_\alpha,u_{\beta \,\alpha})_{\beta \,\alpha \in I  }$ such that $M=\varinjlim F_\alpha$.

\begin{Lemma} \label{zerosum} Let $\B$ be a class of right
$R$-modules closed under direct sums, and let $\Q$ be a class of left $R$-modules.
Assume that   $M$ is
$\B$-stationary.
If for any
pair $\alpha,\beta\in I$ with  $\beta\ge \alpha$ and for any $Q\in
\Q $ there exists $B=B_{\beta \,\alpha}(Q)\in \B $   such that
$$\Ker (u_{\beta\alpha} \otimes_R Q ) =\bigcap _{
h\in H_{u_{\beta\alpha}}(B)} \Ker (h\otimes _RQ), $$ then $M$ is a
$\Q$-Mittag-Leffler module.
\end{Lemma}

\begin{Proof} Let $\Q'=\{Q_k\}_{k\in K}$ be any family of modules in
$\Q$. To prove the statement, we verify that $\Q'$
satisfies the assumption of Lemma~\ref{zero} for \[B=\bigoplus _{Q\in
\Q'}\,\bigoplus_{\alpha, \beta \in I, \beta \ge
\alpha } B_{\beta \,\alpha}(Q)\: \in \B.\]
By hypothesis and by the construction of $B$, if we fix a pair $\alpha,\beta\in I$ with $\beta \ge
\alpha \in I$, then for all $Q\in
\Q' $
$$\Ker (u_{\beta\alpha} \otimes_R Q ) =\bigcap _{
h\in H_{u_{\beta\alpha}}(B)} \Ker (h\otimes _RQ) $$  As $M$ is $B$-stationary, we conclude from Lemma~\ref{zero}
that $M$ is $\Q'$-Mittag-Leffler.
\end{Proof}

\begin{Prop} \label{Prop1} Let $\B$ be a class of right
$R$-modules closed under direct sums, and let $\Q$ be a class of left $R$-modules. Assume that   $M$ is
$\B$-stationary. If for every
$Q\in\Q$ and every $\alpha\in I$ there exists a map
$f_\alpha\colon F_\alpha\to B_\alpha$ such that $B_\alpha\in\B$
and $f_\alpha\otimes_RQ$ is a monomorphism, then $M$ is a
$\Q$-Mittag-Leffler module.
\end{Prop}
\begin{Proof}
We verify the condition in Lemma \ref{zerosum}. Let $\beta\ge
\alpha$ in $I $ and $Q\in\Q$. By hypothesis, there is
$f_\beta\colon F_\beta\to B_\beta\in \B$ such that $f_\beta\otimes_RQ$
is a monomorphism. Set $h_\beta=  f_\beta\, u_{\beta\alpha}$. Then
$h_\beta\in H_{u_{\beta\alpha}}(B_\beta)$, so
$$\bigcap _{
h\in H_{u_{\beta\alpha}}(B_\beta)} \Ker (h\otimes _RQ) \subseteq
\Ker (h_\beta\otimes _RQ)= \Ker (u_{\beta\alpha} \otimes_R Q ) $$
and the reverse inclusion is always true.
\end{Proof}


We have seen several conditions implying that a $\B$-stationary module is $\Q$-Mittag-Leffler. Let us now discuss the reverse implication.
We will need the following notion.

\begin{Def}{\rm Let $\B$ be a class of right $R$-modules, and let $A$ be a right $R$-module.
A morphism $f \in \Hom _R(A,B)$ with $B \in \B$ is a $\B$-{\em
preenvelope} (or a left $\B$-approximation) of $A$ provided that
the abelian group homomorphism $\Hom _R(f, B')\colon \Hom _R(B,B')
\to \Hom _R(A,B')$ is surjective for each $B'\in \B$. }
\end{Def}

\begin{Lemma}\label{E}
Let $\B$ be a class of right $R$-modules, and let $u\colon  M\to N$ and $v\colon  M\to M'$ be right $R$-module
homomorphisms.  Assume that $M'$ has a $\B$-preenvelope $f\colon
M'\to B$. Consider the
following statements.
\begin{itemize}
\item[(1)]  $v$ $\B$-dominates $u$.
\item[(2)] $\Ker (u \otimes_R  B^\ast) \subseteq \Ker (fv\otimes B^\ast) $.
\item[(3)] $\Ker (u \otimes_R B^\ast ) \subseteq \Ker (v\otimes _RB^\ast) $.
\end{itemize}
Statements (1) and (2)  are equivalent, and  statement (3) implies (1) and (2).
Moreover, if  there is a class of left $R$-modules $\Q$ such that
 the character module $B^\ast\in\Q$ and  $f\otimes _RQ$ is a monomorphism for all $Q\in\Q$,
then all three statements  are  equivalent to
\begin{itemize}
\item[(4)] $v$ dominates $u$ with respect to $\Q$.
\end{itemize}
\end{Lemma}
\begin{Proof}
$(1)\Rightarrow (2).$
Since $\widetilde{h}=f\,v\in H_v(B)$, we have $\Ker (u\otimes _RQ) \subseteq \Ker (\widetilde{h}\otimes _RQ)$ for all left $R$-modules $Q$, so in particular for $Q=B^\ast$.

$(2)\Rightarrow (1).$ By \ref{trick} we have $\Ker (u\otimes _RQ) \subseteq \Ker (\widetilde{h}\otimes _RQ)$ for all left $R$-modules $Q$. Let $B'\in\B$.
Since every $h\in H_v(B')$ factors through $\widetilde{h}$, we further have $\Ker (\widetilde{h}\otimes _RQ)\subseteq \Ker ({h}\otimes _RQ)$ for all $h\in H_v(B')$ and all $_RQ$, hence (1) holds true.

$(3)\Rightarrow (2)$ holds true because
$\Ker (v\otimes _RB^\ast) \subseteq \Ker (fv\otimes B^\ast) $.

Assume now that    $B^\ast\in\Q$ and $f\otimes _RQ$ is a monomorphism for all $Q\in\Q$. We prove
$(1)\Rightarrow (4).$ As above we see $\Ker (u\otimes _RQ) \subseteq \Ker (\widetilde{h}\otimes _RQ)$
for all left $R$-modules $Q$. Moreover, if $Q\in\Q$, then $\Ker (\widetilde{h}\otimes _RQ)=\Ker (v\otimes _RQ)  $, which yields (4).

Finally,  $(4)\Rightarrow (3)$ as $B^\ast\in\Q$.
\end{Proof}

\begin{Prop} \label{Prop2}
Let $\B$ be a class of right
$R$-modules.  Assume that  for every $\alpha\in I$ there exists a
$\B$-preenvelope $f_\alpha\colon F_\alpha\to B_\alpha$.
Assume further that
$M$ is  $Q$-Mittag-Leffler for  $Q=\bigoplus _{\alpha\in I}
B_\alpha^\ast$.  Then $M$  is $\B$-stationary.
\end{Prop}
\begin{Proof}
For any $\alpha \in I  $, denote by $u_\alpha \colon  F_\alpha \to
M$ the canonical map. By Theorem \ref{2.1.4} we must show that
 there exists
$\beta \ge \alpha $ such that $u_{\beta \alpha}$
$\B$-dominates $u_{\alpha}$, which  means  $\Ker (u_\alpha\otimes _RB_\beta^\ast)\subseteq \Ker (
u_{\beta \alpha}\otimes _RB_\beta^\ast)$ by Lemma \ref{E}.
So,    it is enough to find $\beta \ge \alpha $ such that
$$\Ker (u_\alpha\otimes _RQ)\subseteq \Ker (u_{\beta \alpha}\otimes _RQ)$$
To this end, we take a generating set $(x_k)_{k\in K}$ of $\Ker (u_\alpha\otimes _RQ)$ and consider the diagram
$$\begin{array}{ccccc}&&F_\alpha\otimes \prod _{k\in K}Q & {\stackrel{u_\alpha\otimes
\prod _{k\in K}Q
}{\longrightarrow}} & M\otimes \prod _{k\in K}Q\\
&&\rho_\alpha\downarrow \phantom{\rho_\alpha}& & \rho  \downarrow \phantom{\rho }\\
&& \prod _{k\in K}(F_\alpha\bigotimes_RQ)&{\stackrel{\prod _{k\in
K}(u_\alpha\otimes Q)}{\longrightarrow}} & \prod _{k\in
K}(M\bigotimes_RQ).\end{array}
$$
Since $\rho_\alpha$ is an isomorphism, there is $x\in F_\alpha\otimes \prod _{k\in K}Q $ such that
$\rho_\alpha (x)=(x_k)_{k\in K}\in \prod _{k\in K}(F_\alpha\bigotimes_RQ)$. Then $(u_\alpha\otimes
\prod _{k\in K}Q) (x)=0$ because $\rho$ is injective, and thus $x\in  \Ker (u_{\beta\alpha}\otimes
_R\prod _{k\in K}Q)$ for some $\beta\ge \alpha$. {}From the diagram
$$\begin{array}{ccccc}&&F_\alpha\otimes \prod _{k\in K}Q & {\stackrel{u_{\beta \alpha}\otimes
\prod _{k\in K}Q
}{\longrightarrow}} & F_\beta\otimes \prod _{k\in K}Q\\
&&\rho_\alpha\downarrow \phantom{\rho_\alpha}& & \rho_\beta  \downarrow \phantom{\rho_\beta }\\
&& \prod _{k\in K}(F_\alpha\bigotimes_RQ)&{\stackrel{\prod _{k\in
K}(u_{\beta \alpha}\otimes Q)}{\longrightarrow}} & \prod _{k\in
K}(F_\beta\bigotimes_RQ).\end{array}
$$
we deduce that $(x_k)_{k\in K}=\rho_\alpha (x)\in \Ker \prod _{k\in K} (u_{\beta\alpha}\otimes _RQ)$, that is,
 $x_k\in  \Ker (u_{\beta\alpha}\otimes _RQ)$ for all ${k\in K}$, and we conclude  $\Ker (u_\alpha\otimes _RQ)\subseteq \Ker (u_{\beta \alpha}\otimes _RQ)$.
\end{Proof}

The previous observations are subsumed in the following result.

\begin{Th} \label{ponte}
Let $\B$ be a class of right
$R$-modules closed under direct sums, and let $\Q$ be a class of left $R$-modules.
Assume that  for every finitely presented module $F$ there exists
a $\B$-preenvelope $f\colon F\to B$ such that the character module
$B^\ast\in\Q$ and  $f\otimes _RQ$ is a monomorphism for all $Q\in\Q$.
Then the following statements are equivalent for a right
$R$-module  $M$.
\begin{itemize}
\item[(1)] $M$ is $\B$-stationary.
\item[(2)] $M$ is  $Q$-Mittag-Leffler  for all $Q\in\text{\rm Add}\Q$.
\item[(3)] $M$ is  $\Q$-Mittag-Leffler.
\end{itemize}
\end{Th}
\begin{Proof}
(1) implies (3) by \ref{Prop1}, (2) implies (1) by \ref{Prop2}, and (3) implies (2) by \ref{definable}.
\end{Proof}

\begin{Ex}\label{endonoeth}
Let $B$ be a right $R$-module with the property that all finite matrix subgroups of $B$ are finitely generated over the endomorphism ring of $B$,
for example an endonoetherian module.
Then a right $R$-module $M$ is  $B^\ast$-Mittag-Leffler if and only if it  is $B$-stationary.
\end{Ex}
\begin{Proof}
Set $\B=\text{Add}\,B$ and $\Q=\text{Add}\,B^\ast$. By Theorem \ref{definable} and Corollary \ref{closureB} we know that
 $M$ is  $B^\ast$-Mittag-Leffler if and only if it  is $\Q$-Mittag-Leffler, and $M$ is $B$-stationary if and only if it  is $\B$-stationary.
Moreover, by \cite[3.1]{via}
   every finitely presented module $F$ has a $\B$-preenvelope $f\colon F\to B'$ with $B'\in\text{add}\,B$, hence $(B')^\ast\in\Q$.
   Finally, $\Hom _R(f,B)$ is an epimorphism, thus applying $\Hom _Z(-,\Q/\Z)$ and using Hom-$\otimes$-adjointness,
   we see that $f\otimes _RQ$ is a monomorphism for all $Q\in\Q$. So the claim follows immediately from \ref{ponte}.
\end{Proof}

Further applications of Theorem \ref{ponte} are given in Section \ref{cotorsion}.

\section{Baer modules}\label{baer}

Throughout this section, let $R$ be a commutative domain. A module $M$ is said to be a {\it Baer
module} if $\mathrm{Ext}_R^1(M,T)=0$ for any torsion $R$-module
$T$.

Kaplansky in 1962 proposed the question whether the only Baer
modules are the projective modules. He was inspired by the
analogous question raised by Baer for the case of abelian groups, which was solved by  Griffith in 1968.

In the general case of domains, a positive answer to Kaplansky's question was recently given by the authors in joint work with S. Bazzoni \cite{abh}.
The proof uses an important   result of Eklof, Fuchs
and Shelah from 1990 which reduces the problem  to showing that countably presented Baer modules are projective (cf.
\cite[Theorem~8.22]{fuchssalce}).

Aim of this section is to
give a new proof for the fact that every countably generated Baer module is projective, which uses our previous results.
In fact, we
 are   going to see that a
countably generated Baer module is Mittag-Leffler. Then the result
follows because countably generated flat Mittag-Leffler modules
are projective (cf.~\cite[Corollaire~2.2.2 p.~74]{RG}).

So, let us consider a countably presented  Baer module $M$. By Kaplansky's work, we know that $M$ is flat of projective dimension at most one. So $M$
can
be written as a direct system of the form
\[F_1\stackrel{f_1}\to F_2\stackrel{f_2}\to F_3\to\dots\to F_{n}\stackrel{f_n}\to F_{n+1}\to\dots\]
where, for each $n\ge1$, $F_n$ is a finitely generated free
$R$-module. As the class of torsion modules is closed under direct
sums, it follows from \cite{bazher}, see \ref{BH},  that $M$ is a Baer module if
and only if for any torsion module $T$ the inverse system
$\left(\mathrm{Hom}_R(F_n,T), \mathrm{Hom}_R(f_n,T)\right)$ satisfies the
Mittag-Leffler condition, in other words, iff $M$ is $T$-stationary.

\begin{Lemma}\label{clau} Let $R$ be a commutative domain.
\begin{itemize}
\item[(i)] Let $Q$ be a finitely generated torsion module. For any $n\ge
0$ there exist a torsion module $T$ and a homomorphism $h\colon
R^n\to T$ such that $h\otimes Q$ is injective.

\item[(ii)] Let $Q$ be a finitely generated torsion-free $R$ module. For
any $n\ge 0$   there exists a torsion module $T$ such that
\[\bigcap_{h\in \mathrm{Hom}(R^n,T)} \Ker(h\otimes Q)=0.\]
\end{itemize}

\end{Lemma}

\begin{Proof} (i). If $n=0$, the claim is trivial (we assume that $0$ is torsion). Fix $n\ge 1$. As $Q$ is finitely
generated and torsion $I=\mathrm{ann} _R(Q)$ is a nonzero ideal of
$R$. So that, $T=(R/I)^n$ is a torsion module. Note that
$R/I\otimes Q\cong Q/IQ=Q$. Hence the canonical projection $R^n\to
T$ satisfies the desired properties.

(ii). First we show that for any $0\neq x\in R^n\otimes Q$ there
exist a torsion module $T_x$ and $h\colon  R^n\to T_x$ such that
$(h\otimes Q)(x)\neq 0$. Let $K$ denote the field of quotients of
$R$. As $Q$ is torsion-free and finitely generated, it can be
identified with a finitely generated submodule of $K^m$ for some
$m$. Moreover, as $Q$ is finitely generated, multiplying   by a
suitable nonzero element of $R$, we can assume that $Q\le R^m$.

The claim is trivial for $n=0$. Fix $n\ge 1$. As $0\neq x\in
R^n\otimes Q=Q^n$ there is $i\in \{1,\dots ,n\}$ such that the
$i$-th component of $x$ is nonzero. Let $\pi _i\colon  R^n\to R$
denote the projection on the i-th component. As $(\pi _i\otimes
Q)(x)\neq 0$, we only need to prove the statement   for $n=1$.

Set $x=(x_1,\dots ,x_m)\in R\otimes Q=Q\subseteq R^m$. Let $j\in
\{1,\dots ,m\}$ be such that $x_j\neq 0$. As $R$ is a domain,
there exists $0\neq t\in R$ such that $tR\subsetneq x_jR$. Hence
$x\not \in tQ\subseteq tR^m$. That is, if $p\colon  R\to R/tR$
denotes the canonical projection, $x$ is  not in $tQ
=\Ker (Q{\stackrel{p\otimes Q }{\longrightarrow}}
R/tR\otimes Q)$.

To prove statement (ii) take  $T=\oplus _{x\in R^n\otimes
Q\setminus \{0\}}T_x$.
\end{Proof}

\begin{Th} If $R$ is a commutative domain then any countably presented  Baer module over $R$ is
Mittag-Leffler. Therefore any Baer module is projective.
\end{Th}

\begin{Proof} Let $M_R$ be a countably presented Baer module. Then
$M_R$ is flat, hence a direct limit of finitely generated free
modules.

Denote by $\Tcal$ and $\Fcal$ the classes of torsion and
torsion-free modules, respectively. By Theorem~\ref{BH},  $M_R$ is
$\mathcal{T}$-stationary. Since $\mathcal{T}$ is closed under direct
sums, the previous Lemma \ref{clau}(i) together with \ref{Prop1}  implies that
$M$ is $\Q$-Mittag-Leffler where $\Q$ is the class of
 all
finitely generated   modules from $\Tcal$.
By Theorem~\ref{definable},  it follows that $M$ is
$\Tcal$-Mittag-Leffler.

Next, we show that $M$ is also $\Fcal$-Mittag-Leffler.
Again by Theorem~\ref{definable},  it is enough to show
that $M$ is
Mittag-Leffler with respect to the class of finitely generated
torsion-free modules. So, let $Q$ be a finitely generated
torsion-free module, and let $u\in\Hom _R(F,F')$ be a morphism between finitely generated free modules $F,F'$.
By Lemma~\ref{clau}(ii)
there  exists a torsion module $T$ such that
$\bigcap_{h\in \mathrm{Hom}(F',T)} \Ker (h\otimes Q)=0.$  So,
 if $x\in F\otimes _RQ$ and $y=(u \otimes_R Q) (x)\not=0$, then there must exist $h'\in\Hom _R(F',T)$
  such that $(h'\otimes _RQ)(y)\not=0$, which means $(h'u\otimes _RQ)(x)\not=0$ and shows that $x\not\in\bigcap _{
h\in H_{u}(T)} \Ker (h\otimes _RQ)$. Thus we deduce that
$\Ker (u \otimes_R Q ) =\bigcap _{
h\in H_{u}(T)} \Ker (h\otimes _RQ)$. Our claim then follows  from Lemma~\ref{zerosum}.

Since $M$ is flat,  we now conclude from Corollary~\ref{tp} that $M$
is Mittag-Leffler and thus projective.
\end{Proof}

\section{Matrix subgroups}\label{smatrix}

In \cite{RG}  Raynaud and Gruson  also studied  modules satisfying a stronger  condition, which they called {\it strict Mittag-Leffler modules}.
In this section, we investigate the relative version of this condition and interpret it in terms of matrix subgroups. Hereby we establish a
relationship with work of W.~Zimmermann \cite{Z}.

We start out with a stronger version of Proposition \ref{dominatinglimits}.

\begin{Prop}\label{strictdominatinglimits}
Let $B$ be a right $R$-module,  and let $\mathcal{S}$ be
a class of finitely  presented right $R$-modules. For a right
$R$-module $M\in \varinjlim \mathcal{S}$, the following statements are equivalent.

\noindent{\textrm{(1).}} There
is a direct system of  finitely presented right $R$-modules
$(F_\alpha,u_{\beta \,\alpha})_{\beta \,\alpha \in I  }$    with
$M=\varinjlim (F_\alpha,u_{\beta \alpha})_{\beta , \alpha \in I }$
having the property that for  any $\alpha \in I  $ there exists
$\beta \ge \alpha$ such that the canonical map $u_{\alpha}:F_{\alpha}\to M$
satisfies  $H_{u_{\alpha}}(B)=H_{u_{\beta \alpha}}(B)$.

\noindent{\textrm{(2).}} Every
direct system of  finitely presented right $R$-modules
$(F_\alpha,u_{\beta \,\alpha})_{\beta \,\alpha \in I  }$    with
$M=\varinjlim (F_\alpha,u_{\beta \alpha})_{\beta , \alpha \in I }$
has the property that for  any $\alpha \in I  $ there exists
$\beta \ge \alpha$ such that the canonical map $u_{\alpha}:F_{\alpha}\to M$ satisfies $H_{u_{\alpha}}(B)=H_{u_{\beta \alpha}}(B)$.

\noindent{\textrm{(3).}} For any finitely presented module $F$ (belonging to $\mathcal S$) and
any homomorphism $u\colon  F\to M$ there exist a   module $S\in
\mathcal{S}$ and a homomorphism $v\colon  F\to S$ such that
 $u$ factors through $v$,
and $H_u(B)=H_v(B)$.

\end{Prop}

\begin{Proof}
Adapt the proof of Proposition \ref{dominatinglimits} replacing Lemma~\ref{composition} by Lemma \ref{diagram}(2) and Lemma \ref{hpure}(1).
\end{Proof}

In view of the characterization of $B$-stationary modules given in Theorem \ref{2.1.4}, we introduce the following terminology.

\begin{Def} {\rm Let $B$ be a right
$R$-module. We say that a right $R$-module $M$ is {\it strict $B$-stationary} if it satisfies the equivalent conditions
of Proposition \ref{strictdominatinglimits} (for    some class of finitely presented modules $\mathcal{S}$).

If $\B$ is a class of right $R$-modules,  then we say that
 $M$ is strict $\B$-stationary if it is strict $B$-stationary  for every $B\in\B$.}
\end{Def}

\begin{Remark}\label{firstremarks} {\rm
(1) The modules that are strict $B$-stationary for every right $R$-module $B$  are exactly the strict
Mittag-Leffler modules of \cite{RG}. In particular, every pure-projective module is strict Mod$R$-stationary, see \cite[2.3.3]{RG}.
\\
(2) By Proposition \ref{2.1.1}  and Theorem \ref{2.1.4}, every strict $B$-stationary module is
$B$-stationary. \\
(3) By Lemma \ref{MLstrict} a countably presented module is strict
$B$-stationary if and only if it is $B$-stationary}.
\end{Remark}

The class $\B$ in the definition of a strict $\B$-stationary module enjoys slightly
weaker closure properties with respect to the class $\B$ in the definition of  $\B$-stationary modules.

\begin{Prop}\label{prodstrictstationary} Let $\{B_j\}_{j\in J}$ be a
family of right $R$-modules. Let $(F_\alpha,u_{\beta
\,\alpha})_{\beta \,\alpha \in I  }$ be a direct system of
finitely presented right $R$-modules  and $M=\varinjlim F_\alpha$.
Then the following statements are equivalent.
\begin{itemize}
\item[(1)] $M$ is strict $\prod _{j\in
J}B_j$-stationary.
\item[(2)] For any $\alpha \in I$ there exists $\beta \ge \alpha$ such
that the canonical map $u_{\alpha}:F_{\alpha}\to M$ satisfies
$H_{u_{\alpha}}(B_j)=H_{u_{\beta \alpha}}(B_j)$ for any $j\in J$.
\item[(3)] $M$ is strict $\bigoplus _{j\in
J}B_j$-stationary.
\end{itemize}
\end{Prop}

\begin{Proof} Proceed as in Proposition \ref{prodstationary}, cf.~Remark \ref{hsubgr}.
\end{Proof}

\begin{Cor}\label{strictclosureB}    Let $\B$ be a class
of right $R$-modules. Let $M$ be a strict $\B$-stationary  right
$R$-module. Then the following statements hold true.
\begin{itemize}
\item[(i)] $M$ is strict $\B '$-stationary where $\B '$ denotes the class
of all modules isomorphic either to a locally split submodule  or to a
pure quotient  of a module  in $\B$.

\item[(ii)] $M$ is strict $\Add \,\B$-stationary if and only if it
is  strict $\Prod\, \B$-stationary if and only if there exists a
direct system of finitely presented right $R$-modules
$(F_{\alpha},u_{\beta \alpha})_{\beta \alpha \in I}$ with
$\varinjlim F_{\alpha}\cong M$ having the property that for any
$\alpha \in I$ there exists $\beta \ge \alpha$ such that the
canonical map $u_{\alpha}:F_{\alpha}\to M$ satisfies
$H_{u_{\alpha}}(B)=H_{u_{\beta \alpha}}(B)$ for any $B\in \B$.

\item[(iii)] $M$ is strict $\Add B$- and strict $\Prod B$-stationary for every $B\in \B$.

\end{itemize}

\end{Cor}
\begin{Proof}
Adapt the proof of Corollary \ref{closureB}   replacing
 Lemma~\ref{hpure}(2)    by
  Lemma~\ref{hpure}(4), and     Proposition~\ref{prodstationary} by  Proposition~\ref{prodstrictstationary}.
  \end{Proof}

We now recall some notions from \cite{Z}.

\begin{Def}\label{defmatrix} {\rm Given two right $R$-modules $A,B$, an integer $n\in\N$, and an element $\underline{a}=(a_1,\ldots,a_n)\in A^n$,
we consider the $\End B$-linear map
$$\varepsilon_{\underline{a}}: \Hom _R(A,B) \to B^n,\, f\mapsto f(\underline{a})=(f(a_1),\ldots,f(a_n))$$
and define $$H_{{A,\underline{a}}}(B) = \text{Im} \,\varepsilon_{\underline{a}}= \{f(\underline{a})\,\mid\,f\in\Hom _R(A,B)\}.$$

\noindent
If $n=1$, then $H_{{A,\underline{a}}}(B)=H_{{A,{a}}}(B)$ is called a {\it matrix subgroup} of $B$, and it is called  a
{\it finite matrix subgroup} if the module $A$ is finitely presented.
 }\end{Def}

 The subgroups $H_{{A,\underline{a}}}(B)$ are related to $H$-subgroups as follows.

 \begin{Lemma}\label{Hmatrix}
Let $A,B,M$ be  right $R$-modules,   $n\in\N$, and   $\underline{a}=(a_1,\ldots,a_n)\in A^n$. \\
(1) If $a_1,\ldots,a_n$ is a generating set of $A$, then $\varepsilon_{\underline{a}}$ is a monomorphism.\\
(2) If $v\in\Hom _R(A,M)$, then $\varepsilon_{\underline{a}}(H_v(B))=H_{{M,\underline{m}}}(B) $ where
$\underline{m}=v(\underline{a})\in M^n$.\\
(3) If $\underline{m}=(m_1,\ldots,m_n)\in M^n$, then $H_{{M,\underline{m}}}(B)=\varepsilon_{\underline{e}}(H_u(B)) $ where
$\underline{e}=(e_1,\ldots,e_n)$ is given by the canonical basis $e_1,\ldots,e_n$ of the free module $R^n$,
and $u:R^n\to M$ is defined by $u(r_1,\ldots,r_n)=\sum_{i=1}^n m_i\,r_i$.
 \end{Lemma}
 \begin{Proof} Is left to the reader.\end{Proof}

 We will need some further terminology from \cite{Z}.

 \begin{Def} {\rm Let  $n\in\N$. A pair $(A,\underline{a})$ consisting of  a right $R$-module $A$ and an element $\underline{a}=(a_1,\ldots,a_n)\in A^n$
 will be called an \emph{$n$-pointed module}. A \emph{morphism of $n$-pointed modules} $h:(A,\underline{a})\to (M,\underline{m})$ is an $R$-module
 homomorphism $h:A\to M$ such that $h(\underline{a})=\underline{m}$.

Consider  now
a direct system of   right $R$-modules $(F_\alpha,u_{\beta \,\alpha})_{\beta
\,\alpha \in I  }$    with direct limit $M=\varinjlim (F_\alpha,u_{\beta \alpha})_{\beta , \alpha \in I }$ and canonical maps $u_\alpha:F_\alpha\to M$.
If for every $\alpha\in I$ the elements $\underline{x}_\alpha\in F_\alpha\,^n$ are chosen in such a way that the
$u_{\beta \alpha}: (F_\alpha,\underline{x}_\alpha)\to (F_\beta,\underline{x}_\beta)$
are morphisms of $n$-pointed modules, then $( (F_\alpha,\underline{x}_\alpha), u_{\beta \alpha})_{\beta,\alpha\in I}$ is called a \emph{direct system of $n$-pointed modules}. Setting $\underline{m}=u_\alpha(x_\alpha)$ for some
 $\alpha\in I$, we have that also the $u_\alpha:(F_\alpha, \underline{x}_\alpha)\to (M, \underline{m})$ are morphisms of $n$-pointed modules.
 We then write $ (M, \underline{m})=\varinjlim (F_\alpha, \underline{x}_\alpha)$.
}\end{Def}

We now show that the strict $B$-stationary modules are precisely the modules studied by Zimmermann
in \cite[3.2]{Z}. Like the Mittag-Leffler modules, they can be characterized in terms of the
injectivity of a natural transformation.

Let ${}_SB_R$ be and $S$-$R$-bimodule, and let ${}_SV$ be a left $S$-module. For any right
$R$-module $M_R$ there is a natural transformation
\[\nu =\nu {(M,B,V)} \colon M\otimes _R\mathrm{Hom}_S(B,V)\to \mathrm{Hom}_S(\mathrm{Hom}_R(M,B),V)\]
defined by
$$\nu(m\otimes \varphi):\; f\mapsto \varphi (f(m)).$$
Notice that when $B=V$ and $_SB_R$ is faithfully balanced, then $\nu$ is induced by the evaluation map of $M$ inside its bidual.

If $M_R$ is finitely presented and ${}_SV$ is injective then $\nu$ is an isomorphism
(cf.~\cite[Theorem~3.2.11]{enochsjenda}). The case that  $\nu$ is a monomorphism for all injective modules $_SV$ was studied by Zimmermann in
\cite[3.2]{Z}. We are going to see below that this happens precisely when $M$ is  strict
$B$-stationary. So, like   $Q$-Mittag-Leffler modules, strict $B$-stationary
modules can be characterized in terms of the injectivity of a natural transformation, in this case
the injectivity of $\nu$. Let us first discuss how the injectivity of $\nu$ behaves
under direct sums.

If $M=\oplus _{i\in I}M_i$ then, for each $i\in I$, the canonical inclusion $M_i\to M$ induces an
inclusion
\[\mathrm{Hom}_S(\mathrm{Hom}_R(M_i,B),V) \to \mathrm{Hom}_S(\mathrm{Hom}_R(M,B),V).\]
This family of inclusions induces an injective map
\[\Phi: \oplus _{i\in I}\mathrm{Hom}_S(\mathrm{Hom}_R(M_i,B),V) \to
\mathrm{Hom}_S(\mathrm{Hom}_R(M,B),V)\]
given by the rule
$$\Phi((g_i)_{i\in{I_0}}):\;f\mapsto\sum_{i\in{I_0}} g_i(f\!\!\mid_{M_i}).$$
This allows us to deduce that

\begin{Lemma}\label{nuds} The map $\nu (\oplus _{i\in I}M_i,B,V)$ is injective if and only if $\nu
(M_i,B,V)$ is injective  for any $i\in I$.
\end{Lemma}

 \begin{Th}\label{dcc}
 Let $B$ and $M$ be right $R$-modules. The following statements are equivalent.
 \begin{itemize}
\item[(1)] $M$ is strict $B$-stationary.

\item[(2)] For every $n\in\N$, every element $\underline{m}\in M^n$ and every direct system of $n$-pointed modules
$( (F_\alpha,\underline{x}_\alpha), u_{\beta \alpha})_{\beta,\alpha\in I}$ with all $F_\alpha$ being finitely presented
 and  $ (M, \underline{m})=\varinjlim (F_\alpha, \underline{x}_\alpha)$
there is $\beta\in I$ such that $H_{{M,\underline{m}}}(B)=H_{{F_\beta,\underline{x}_\beta}}(B)$.

\item[(3)] For every $n\in\N$ and every element $\underline{m}\in M^n$ there are a finitely presented right $R$-module $A$, an element
 $\underline{a}\in A^n$ and a morphism of $n$-pointed modules  $h:(A,\underline{a})\to (M,\underline{m})$ such that
 $H_{{M,\underline{m}}}(B)=H_{{A,\underline{a}}}(B)$.
 \end{itemize}

 \noindent
 Let $S$ be a ring such that $_SB_R$ is a bimodule. Then the following statement is further equivalent.
 \begin{itemize}

\item[(4)] For every injective left $S$-module $V$, the canonical map $\nu: M\otimes_R\, \Hom _S (B,V)\to  \Hom _S (\Hom _R (M,B), V)$
defined by $\nu(m\otimes \varphi): f\mapsto \varphi (f(m))$
 is a monomorphism.
 \end{itemize}
 \end{Th}
 \begin{Proof}
 The equivalence of the last three statements is due to Zimmermann, see \cite[3.2]{Z}.

 (1) $\Rightarrow$ (3):
 As in Lemma \ref{Hmatrix}(3), we write
 $H_{{M,\underline{m}}}(B)=\varepsilon_{\underline{e}}(H_u(B)) $
 where
$\underline{e}=(e_1,\ldots,e_n)$ is given by the canonical basis $e_1,\ldots,e_n$ of the free module $R^n$,
and $u:R^n\to M,\, (r_1,\ldots,r_n)\mapsto\sum_{i=1}^n m_i\,r_i$.
By assumption there exist a  finitely presented module $A$, a homomorphism $v\colon  R^n\to A$ and a homomorphism $h: A\to M$ such that
 $u=h\,v$,
and $H_u(B)=H_v(B)$.
 Set $\underline{a}=v(\underline{e})\in A^n$. Then
 $h(\underline{a})= u(\underline{e})=\underline{m}$. So, we obtain a morphism of $n$-pointed modules $h: (A,\underline{a})\to (M,\underline{m})$. Moreover, we see as in Lemma \ref{Hmatrix}(2) that $H_{{M,\underline{m}}}(B)=\varepsilon_{\underline{e}}(H_v(B))= H_{{A,\underline{a}}}(B)$.

 (2) $\Rightarrow$ (1):
 Take a direct system of  finitely presented right $R$-modules
$(F_\alpha,u_{\beta \,\alpha})_{\beta \,\alpha \in I  }$    with
$M=\varinjlim (F_\alpha,u_{\beta \alpha})_{\beta , \alpha \in I  }$\,, and fix $\alpha \in I  $.
Moreover, choose a generating set $a_1,\ldots,a_n$ of $F_\alpha$, and set $\underline{x}_\alpha=(a_1,\ldots,a_n)$. Set further
$\underline{x}_\beta=u_{\beta\alpha}(\underline{x}_\alpha)$ for all $\beta\ge \alpha$, and $\underline{m}=u_{\alpha}(\underline{x}_\alpha)$. Then $(F_\beta,\underline{x}_\beta)_{\beta\in I, \beta\ge \alpha}$
is a direct system of $n$-pointed modules with direct limit
$(M,\underline{m})$.
So, by assumption there is $\beta\in I$ such that $H_{{M,\underline{m}}}(B)=H_{{F_\beta,\underline{x}_\beta}}(B)$.
By Lemma \ref{Hmatrix}(2) this means
$\varepsilon_{\underline{x}_\alpha}(H_{u_{\alpha}}(B))= \varepsilon_{\underline{x}_\alpha}(H_{u_{\beta\alpha}}(B))$.
Since $\varepsilon_{\underline{x}_\alpha}$ is a monomorphism, we conclude
$H_{u_{\alpha}}(B)=H_{u_{\beta\alpha}}(B)$.
 \end{Proof}

Here are some consequences of the previous theorem.

\begin{Cor}\label{sigmapi}
Let $B$ be a right $R$-module.\\
(1)   \cite[3.6]{Z} The class of
strict $B$-stationary modules is closed under pure submodules and pure extensions. \\
(2) A direct sum of
modules is strict $B$-stationary if and only if so are all direct summands.\\
(3) \cite[3.8]{Z} The module
$B$ is $\Sigma$-pure-injective if and only if every right $R$-module is strict $B$-stationary.
\end{Cor}

The characterization of strict $B$-stationary modules in terms of the injectivity of $\nu$ allows
us to provide a wide class of  examples of such modules. It is the analog of the class discussed in Proposition \ref{MLfilters}.

\begin{Prop}\label{strictfilters} Let $\mathcal{S}$ be a class of right $R$ modules that is strict $\mathcal{B}$-stationary with respect
to a class $\mathcal{B}\subseteq \mathcal{S}^\perp$. Then any
module isomorphic to direct summand of an $\mathcal{S}\cup{\rm
Add}\,R $-filtered module is strict $\mathcal{B}$-stationary.
\end{Prop}

\begin{Proof} As projective modules are strict Mod$R$-stationary and
$\mathcal{S}^\perp=(\mathcal{S}\cup {\rm Add}\,R )^\perp$, we can assume that
$\mathcal{S}$ contains ${\rm Add}\,R$ . By Corollary~\ref{sigmapi}, the class of strict
$\mathcal{B}$-stationary modules is closed by direct summands. So, we only need to prove the statement
for $\mathcal{S}$-filtered modules. Also by Corollary~\ref{sigmapi}, we know that arbitrary direct sums of
modules in $\mathcal{S}$ are strict $\B$-stationary.

Let $M$ be an $\mathcal{S}$-filtered right $R$-module. Let $\tau$ be an ordinal such that there
exists an $\mathcal{S}$-filtration $(M_\alpha)_{\alpha \leq \tau}$ of $M$. Observe that for any
$\beta \leq \alpha \leq \tau$, $M_\alpha$ and $M_\alpha /M_\beta$ are $\mathcal{S}$-filtered
modules, so they belong to ${}^\perp \mathcal{B}$ by \cite[3.1.2]{GT}. For the rest of the proof we fix $B\in
\mathcal{B}$, a ring $S$ such that ${}_SB_R$ is a bimodule, and an injective left $S$-module $V$.
 We shall show that $M$ is strict $\mathcal{B}$-stationary proving  by induction  that
  for any $\alpha \le \tau$ the canonical map
  \[\nu _\alpha \colon M_\alpha \otimes _R\mathrm{Hom}_S(B,V)\to \mathrm{Hom}_S(\mathrm{Hom}_R(M_\alpha,B),V)\]
  is injective.

  As $M_0=0$ the claim is true for
$\alpha =0$. If  $\alpha <\tau$ then the exact sequence
\[0\to M_\alpha \to M_{\alpha +1}\to M_{\alpha +1}/M_\alpha \to 0\]
and the fact that $\mathcal{B}\subseteq \mathcal{S}^\perp$ yields a commutative diagram with exact
rows
$$\begin{array}{ccccccc}&\!M_\alpha \otimes _R\mathrm{Hom}_S(B,V) \! \! & \! \! \to & M_{\alpha +1}\otimes _R\mathrm{Hom}_S(B,V)\! \!&\! \!
{\to}\! \!&\! \!(M_{\alpha +1}/M_\alpha)\otimes _R\mathrm{Hom}_S(B,V)&\!\! \! \to 0\\ &\!\nu
_\alpha\downarrow \phantom{\nu _\alpha}& & \nu _{\alpha +1}  \downarrow \phantom{\nu _{\alpha +1}
}&&
\nu  \downarrow \phantom{\nu }&\!\\
0\to &\! \mathrm{Hom}_S(\mathrm{Hom}_R(M_\alpha,B),V)\! \!&\! \!\to\! \! &\! \!
\mathrm{Hom}_S(\mathrm{Hom}_R(M_{\alpha +1},B),V)\! \!& \! \!{\to}\! \! &\! \!
\mathrm{Hom}_S(\mathrm{Hom}_R(M_{\alpha +1}/M_\alpha,B),V)\! \!&\!\! \! \to 0.\end{array}
$$
The natural map $\nu$ is injective because  $M_{\alpha +1}/M_\alpha$ is strict $\B$-stationary. So, if $\nu _\alpha$ is injective,  then $\nu _{\alpha +1}$
is also injective.

Let $\alpha \leq \tau$ be a limit ordinal, and assume that $\nu _\beta$ is injective for any
$\beta <\alpha$. We shall prove that $\nu _\alpha $ is injective. Let $x\in \mathrm{Ker}\, \nu
_\alpha $. There exists $\beta <\alpha $ and $y\in M_\beta \otimes _R \mathrm{Hom}_S(B,V)$ such
that $x=(\epsilon _\beta \otimes _R \mathrm{Hom}_S(B,V))\, (y)$, where $\epsilon _\beta \colon M_\beta
\to M_\alpha$ denotes the canonical inclusion. As $M_\alpha /M_\beta \in {}^\perp \B$, $\epsilon
_\beta$ induces an injective map $\epsilon \colon \mathrm{Hom}_S(\mathrm{Hom}_R(M_\beta ,B),V) \to
\mathrm{Hom}_S(\mathrm{Hom}_R(M_\alpha ,B),V)$. Considering the commutative diagram
$$\begin{array}{ccccc}&&M_\beta \otimes _R \mathrm{Hom}_S(B,V) & {\stackrel{\epsilon_\beta \otimes _R \mathrm{Hom}_S(B,V)}
{\longrightarrow}} & M_\alpha \otimes _R\mathrm{Hom}_S(B,V)\\
&&\nu _\beta\downarrow \phantom{\nu _\beta}& & \nu _\alpha  \downarrow \phantom{\nu _\alpha }\\
&& \mathrm{Hom}_S(\mathrm{Hom}_R(M_\beta ,B),V)&{\stackrel{\epsilon}{\longrightarrow}} &
\mathrm{Hom}_S(\mathrm{Hom}_R(M_\alpha ,B),V)\end{array}
$$
we see that $0=\epsilon \nu _\beta (y)$. As $\nu _\beta $  and $\epsilon$ are injective,  $y=0$.
Therefore $x=0$, and $\nu _\alpha$ is injective.
\end{Proof}

Next, we investigate the relationship between relative Mittag-Leffler modules and strict stationary modules. It was shown by Azumaya \cite[Proposition 8]{Az3}
that a module $M$ is strict Mittag-Leffler if and only if every pure-epimorphism $X\to M$, $X\in \Mod R$, is locally split.
We will now see that also the dual property plays an important role in this context. According to \cite{Zlpi}, we will say that a right $R$-module $B$
is \emph{locally pure-injective} if every pure-monomorphism $B\to X$,  $X\in \Mod R$,
is locally split.

Moreover, in the following, for a right $R$-module $B$,  we will indicate by $B^{\bullet}$   a left $R$-module which is obtained from $B$ by some duality,
that is, by taking  a ring $S$ such that   $_SB_R$ is a bimodule together with  an injective cogenerator $_SV$ of $S {\LMod}$,
and setting  $_RB^{\bullet} = \Hom _S(B,V)$.
For example, $B^{\bullet}$ can be the character module $B^\ast$ of $B$. But it can also be the {\em local dual} $B^+$  of $B$, which is obtained as above by choosing $S=\End_RB$.
For a left $R$-module $C$,  the notation $C^{\bullet}$ is used correspondingly.

\begin{Prop}\label{dualizing} Let $M$ and $B$ be  right $R$-modules, and let $C$ be a left $R$-module.
The following statements hold true.
 \begin{itemize}
\item[(1)]  \cite[3.3(1)]{Z}  $M$ is a $C$-Mittag-Leffler module if and only if
$M$ is strict $C^{\bullet}$-stationary.
\item[(2)]  If $M$ is strict $B$-stationary, then $M$ is $B^{\bullet}$-Mittag-Leffler. The converse holds true if $B$ is locally pure-injective.
\end{itemize}
\end{Prop}
\begin{Proof}
(2) The first part of the statement is shown by Zimmermann \cite[3.3(2)(a)]{Z}.
For   the converse, assume that $B$ is locally pure-injective.
Let $_RB^{\bullet} = \Hom _S(B,V)$ where $S$ is a ring such that   $_SB_R$ is a bimodule and  $_SV$ is an injective cogenerator of $S {\LMod}$.
Consider   a ring $T$ such that  $_SV_T$ is a bimodule,
 let  $U_T$  be an injective cogenerator
of $\Mod T$, and assume w.l.o.g. that $V_T\subseteq U_T$. Then $_RB^{\bullet}\,_T$ is also a bimodule, and we can consider
 $B^{\bullet}\,^{\bullet}=\Hom _T(B^{\bullet},U)$.
By (1), $M$ is strict $B^{\bullet}\,^{\bullet}$-stationary.
Furthermore, the evaluation map
$B\to B^{\bullet}\,^{\bullet}$ is a pure monomorphism (see e.g. \cite[1.2(4)]{zbcn}), hence locally split.
By Corollary \ref{strictclosureB}(i) it follows that
$M$ is strict $B$-stationary.
\end{Proof}

\begin{Ex}\label{followsstrict}
Let $B$ be a locally pure-injective
  right $R$-module with the property that all finite matrix subgroups of $B$ are finitely generated over the endomorphism ring of $B$.
  Then a right $R$-module $M$ is strict $B$-stationary if and only if it is $B$-stationary.
  \\
  In particular, this applies to the case when  $B$ is a pure-projective right $R$-module over a
  left pure-semisimple ring $R$.
    \end{Ex}
\begin{Proof}
The statement follows by combining Example \ref{endonoeth} with Proposition \ref{dualizing}(2). \\
When
$R$ is a left pure-semisimple ring, all finitely presented right $R$-modules are endofinite \cite{He}. Hence  every pure-projective right $R$-module
    $B$ is  locally pure-injective by \cite[2.4]{Zlpi}, and   endonoetherian  by \cite{ZHZ}. So, the assumptions are satisfied   in this case.
\end{Proof}

Restricting to local duals, we can employ recent work of Dung and Garcia \cite{DG} to obtain a criterion for endofiniteness of finitely presented modules.

\begin{Prop}\label{endofinite}
Let  $_RC$ be a finitely presented  left $R$-module.
The following statements are equivalent.
\begin{itemize}
\item[(1)] $C$ is endofinite.
\item[(2)] $C$ is $\Sigma$-pure-injective, and $C^+$ is $C$-Mittag-Leffler.
\item[(3)]  $C$ is $\Sigma$-pure-injective, and all cyclic $\End C^+$-submodules of $C^+$ are finite matrix subgroups.
\end{itemize}
\end{Prop}
\begin{Proof}
(1)$\Rightarrow$ (2): $C$ is endofinite if and only if it satisfies the descending and the ascending chain condition on finite matrix subgroups. Hence $C$ is $\Sigma$-pure-injective, and   every right $R$-module is $C$-Mittag-Leffler, see Example \ref{extp}(3).
\\
(2)$\Rightarrow$ (3): By Proposition \ref{dualizing} the module
$C^+$ is strict $C^+$-stationary. By condition (3) in Theorem
\ref{dcc} for $n=1$, it follows that all matrix subgroups of $C^+$
of the form $H_{{C^+},\,m}(C^+)$ with $m\in C^+$ are finite matrix
subgroups, and
 of course, the  matrix subgroups of such form are precisely the cyclic $\End C^+$-submodules of $C^+$.
\\
(3)$\Rightarrow$ (1):
Since $C$ is $\Sigma$-pure-injective, the module $C^+$ satisfies the ascending chain condition on finite matrix subgroups, see \cite[Proposition 3]{ZHZ}. Furthermore, every finitely generated $\End C^+$-submodule of $C^+$
is a finite sum of cyclic submodules, hence a finite matrix subgroup, because the class of finite matrix subgroups is closed under finite sums,
see e.g. \cite[2.5]{zbcn}. So, we conclude that  $C^+$ satisfies the ascending chain condition on finitely generated $\End C^+$-submodules,
in other words,  $C^+$ is endonoetherian, see also \cite{Zext}.
Now the claim follows from \cite[4.2]{DG}, where it is shown that a finitely presented module is endofinite provided its local dual is endonoetherian.
\end{Proof}

Using \cite[4.1]{DG1}, we obtain the following observation.

\begin{Cor}\label{criterion}
A left pure-semisimple ring has finite representation type if and only if the local dual $C^+$ of any finitely presented left $R$-module $C$ is $C$-Mittag-Leffler.
\end{Cor}

Before turning in more detail to   pure-semisimple rings, let us apply our results to  the following setting.

\section{Cotorsion pairs}\label{cotorsion}

In this section, we shall see that the theory of relative Mittag-Leffler modules and (strict)
 stationary modules fits very well into the theory of cotorsion pairs.

\begin{Def}\label{cotorsionpair}
{\rm
(1) Let $\Mcal,\Lcal \subseteq \Mod R$
be classes of modules.
The pair  $(\Mcal,\Lcal)$ is said to be a {\em cotorsion pair}
provided
$\Mcal = {}^{\perp} \Lcal$ and $\Lcal = \Mcal ^{\perp}$.

The cotorsion pair $(\Mcal,\Lcal)$ is said to be \emph{complete} if
for every module $X$ there are short exact sequences $0\to X\to L\to
M\to 0$ and $0\to L'\to M'\to X$ where $L,L'\in\Lcal$ and
$M,M'\in\Mcal$

(2) If $\Scal$ is a set of right $R$-modules, we obtain    a
cotorsion pair $(\Mcal,\Lcal)$ by setting $\mathcal L = \Scal
^{\perp}$ and  $\Mcal = {}^{\perp }(\Scal ^{\perp })$. It is called
the cotorsion pair {\em generated}\footnote{This terminology differs
from  previous use, cf.\cite{GT}.}  by $\Scal$, and it is a {\em
complete} cotorsion pair (cf.~\cite[Theorem~3.2.1]{GT}).

(3)
We will say that a cotorsion pair $(\Mcal,\Lcal)$ is
of {\em finite type}
provided it is  generated by a set  of  modules $\Scal\subseteq\rmod R$.
Note that we can always assume $\Scal=\Mcal\cap\rmod R$.

(4) Dually, if $\Scal$ is a set of right $R$-modules, we obtain    a cotorsion pair $(\Mcal,\Lcal)$ by setting
$\Mcal = {}^{\perp }\Scal $ and $\mathcal L = (^{\perp }\Scal) ^{\perp}$. It is called the cotorsion pair {\em cogenerated} by $\Scal$}.
\end{Def}

For more information on cotorsion pairs, we refer to \cite{GT}.

\medskip

Certain classes of complete cotorsion pairs provide a good setting for relative stationarity and Mittag-Leffler properties. As a first
approach we give the following result.

\begin{Prop} \label{ctfirst} Let $(\Mcal,\Lcal)$ be a cotorsion pair
in $\Mod R$. Set $\C=\Mcal^\intercal$. Then the following hold true.
\begin{itemize}
\item[(1)] If $(\Mcal,\Lcal)$ is complete and    $\Lcal$ is closed by direct sums, then any $\Lcal$-stationary right $R$-module is $\C$-Mittag-Leffler.
\item[(2)] If $(\Mcal,\Lcal)$ is generated by (a set of) countably presented modules  and $\Lcal$ is closed by direct sums, then  any module in
$\Mcal$ is strict $\Lcal$-stationary (and thus
$\C$-Mittag-Leffler).
\item[(3)] Assume that $(\Mcal,\Lcal)$ is
generated by a class $\mathcal{S}$ of finitely presented
modules with the property that the first syzygy of any module in $\mathcal{S}$ is also finitely presented.
Then a countably generated module $M$ belongs to $\Mcal$
if and only if
 $M$ belongs to $\varinjlim \mathcal{S}$ and is
(strict) $\Lcal$-stationary.
\end{itemize}
\end{Prop}

\begin{Proof} The hypotheses in $(1)$ imply that any right $R$-module $X$ fits into an exact
sequence  $0\to X\to L\to M\to 0$  where $L\in\Lcal$ and
$M\in\Mcal$.  Hence the statement   follows from
Proposition~\ref{Prop1}.

To prove $(2)$ observe first that a countably presented module in
$\Mcal$ is strict $\Lcal$-stationary by Example \ref{exstat}(1) and Remark \ref{firstremarks} (3).

Now if $M\in \mathcal{M}$ then, by \cite[Corollary~3.2.4]{GT}, $M$
is a direct summand of  a module $N$ filtered by countably presented
modules in $\Mcal$. By the observation above, $N$ is   filtered by strict
$\Lcal$-stationary modules. Then $M$ is strict $\Lcal$-stationary by
Proposition~\ref{strictfilters}. As the cotorsion pair is complete by
\ref{cotorsionpair}(2), we infer from $(1)$  that $M$ is also
$\C$-Mittag-Leffler.

$(3)$
By \cite[2.3]{at2}  $\mathcal{M}$ is included in  $\varinjlim
\Scal$, and, by \cite[Lemma~10.2.4]{enochsjenda}, $\Lcal$ is closed
under direct sums. So, the only-if part follows from (2).
For the converse implication, let $M$ be a countably
generated module in $\varinjlim \Scal$ which is $\Lcal$-stationary.
Statement $(1)$ implies that $M$ is $\Ccal$-Mittag-Leffler. As $R\in
\Ccal$, we deduce from Corollary~\ref{MLcountable} that $M$ is
countably presented. Therefore, and because $\Lcal$ is closed under
direct sums, we can apply Theorem~\ref{BH} to conclude that
$\mathrm{Ext}_R^1(M,L)=0$ for any $L\in \Lcal$. Thus $M\in \Mcal$.
\end{Proof}

>From \cite[1.9]{SaT} we immediately obtain the following consequence.

\begin{Ex}
Let $R$ be an $\aleph_0$-noetherian hereditary ring. If $(\Mcal,\Lcal)$ is a cotorsion pair
in $\Mod R$ such that $\Lcal$ is closed by direct sums, then  any module in
$\Mcal$ is strict $\Lcal$-stationary (and thus
$\C$-Mittag-Leffler).
\end{Ex}

In the following results, we use again the notation $B^{\bullet}$ to
indicate a module obtained from $B$ by some duality, like the
character module, or the local dual of $B$. For a class of modules
$\Scal$, we write $\Scal^{\bullet}$ in order to indicate a class
consisting of  modules that are obtained by some duality from the
modules of $\Scal$. Note that we are not assuming a functorial
relationship between $\Scal$ and $\Scal^{\bullet}$.

\begin{Lemma}\label{ft}
 Let $(\Mcal,\Lcal)$ be a cotorsion pair
of {finite type}
  in $\Mod R$. Set
    $\Scal=\Mcal\cap\rmod R$, and   $\C=\Mcal^\intercal$.
Then the following hold true.
\begin{itemize}
\item[(1)]  $\C=\Scal^\intercal$, and $^\intercal\C=\varinjlim\Scal=\varinjlim\Mcal$.
\item[(2)] If $\D=\C^\perp$, then   $(\C, \D)$ is the cotorsion pair cogenerated by $\Scal^{\bullet}$.
\item[(3)] A right $R$-module $B$ belongs to $\Lcal$ if and only if $B^{\bullet}$ belongs to $\C$.
\item[(4)] A left $R$-module $C$ belongs to $\C$ if and only if $C^{\bullet}$ belongs to $\Lcal$.
\item[(5)] A right $R$-module $B$ belongs to $\varinjlim\Scal$ if and only if $B^{\bullet}$ belongs to $\D$.
\item[(6)] Assume that the cotorsion pair $(\C, \D)$ is of finite type. Then a left $R$-module $X$ belongs to $\D$ iff $X^{\bullet}$ belongs to $\varinjlim\Scal$.
\item[(7)] If $\Ecal=(\varinjlim\Scal)^\perp$, then $(\varinjlim\Scal,\Ecal)$ is the
cotorsion pair cogenerated by the pure-injective modules from $\Lcal$.
\item[(8)] Let $f: N\to M$ be a monomorphism with $M\in\varinjlim\Scal$. Then $f\otimes C$ is a monomorphism for all $C\in \C$
 if and only if $\coker f\in\varinjlim\Scal$.
\end{itemize}
\end{Lemma}
\begin{Proof}
(1) By \cite[2.3]{at2} we have $\Scal\subseteq\Mcal\subseteq\varinjlim\Scal={}^\intercal(\Scal^\intercal)$, hence
$\C=\Mcal^\intercal=\Scal^\intercal$, and $^\intercal\C=\varinjlim\Scal$.
By the well-known Ext-Tor relations we further obtain $\C=\,{^\perp(\Scal^{\bullet})}$, hence (2),  and also statements (3) - (5).
\\
For statement (6), we assume that the cotorsion pair $(\C, \D)$  is of finite type.
Then $X\in\D$ iff $\Ext (C,X)=0$ for all $C\in\mbox{$\C\cap\!\!\mod R$}$, which is equivalent to $\Tor (X^\bullet, C)=0$ for all
$C\in\mbox{$\C\cap\!\!\mod R$}$. But since $\C\subseteq \varinjlim \mbox{$(\C\cap\!\!\mod R)$}$ by \cite[2.3]{at2},
and Tor commutes with direct limits, the latter means that $X^\bullet\in{^\intercal\C}=\varinjlim\Scal$.
\\
Statement (7) is \cite[2.4]{at2}.\\
(8) If $0\to N\mapr{f} M\to Z\to 0$ is exact, then $\Tor_R(M,C)=0$  for all $C\in \C$ by (1). Hence  $f\otimes C$ is a monomorphism for all $C\in \C$
 if and only if $Z\in{}^\intercal \C=\varinjlim \Scal$.
\end{Proof}

As an application of our previous results, we obtain

\begin{Th} \label{cotpair}  Let $(\Mcal,\Lcal)$ be a cotorsion pair
of {finite type}
  in $\Mod R$. Set $\Scal=\Mcal\cap\rmod R$ and
      $\C=\Mcal^\intercal$, and denote by $\Lcal'$   the class of all locally pure-injective modules from $\Lcal$.
Then the following statements are equivalent for a right $R$-module  $M$.
\begin{itemize}
\item[(1)] $M$ is $\Lcal$-stationary.
\item[(2)] $M$ is  $C$-Mittag-Leffler  for all $C\in\C$.
\item[(3)] $M$ is  $\C$-Mittag-Leffler.
\item[(4)] $M$ is strict $C^\bullet$-stationary for all $C\in\C$.
\item[(5)] $M$ is strict $\Lcal'$-stationary.
\end{itemize}
Moreover, every $M\in\Mcal$ is strict $\Lcal$-stationary (and thus $\C$-Mittag-Leffler).  If $M$ is
countably generated, then $M\in\Mcal$ if and only if $M$ belongs to $\varinjlim\Scal$ and   is (strict)
$\Lcal$-stationary.
\end{Th}
\begin{Proof}
First of all, as $\Lcal=\Scal ^\perp$,  the class  $\Lcal$ is closed
under direct sums \cite[Lemma~10.2.4]{enochsjenda}. Moreover, the
cotorsion pair  is complete, see \ref{cotorsionpair}(2). Therefore
for every right $R$-module $F$ there is a short exact sequence $0\to
F\mapr{f} B\to B/F\to 0$ where $B\in\Lcal$ and $B/F\in\Mcal$. Of
course, $f$ is an $\Lcal$-preenvelope. Further, by Lemma \ref{ft},
$B^\ast\in\C$,  and $B/F\in {}^\intercal\C$, hence the map $f\otimes
C$ is a monomorphism for all $C\in\C$. The equivalence of the first
three statements then follows directly from Theorem \ref{ponte}.
Furthermore, (2) is equivalent to (4) by Proposition
\ref{dualizing}.
\\
(2)$\Rightarrow$(5): Let $B\in\Lcal$ be locally pure-injective. Then  $B^{\bullet}\in\C$ by Lemma \ref{ft}, so $M$ is $B^{\bullet}$-Mittag-Leffler, hence strict $B$-stationary by Lemma \ref{dualizing}(2).
\\
(5) $\Rightarrow$(4) is clear since $C^\bullet$ is a (locally) pure-injective module in $\Lcal$ by Lemma \ref{ft} and \cite[1.6(2)]{zbcn}.

The rest of the theorem  follows from Proposition~\ref{ctfirst}.
\end{Proof}

\begin{Cor}\label{CML}
Let $(\Mcal,\Lcal)$ be a cotorsion pair
of {finite type}
  in $\Mod R$. Set
    $\Scal=\Mcal\cap\rmod R$, and   $\C=\Mcal^\intercal$.
    Then a module $M\in\varinjlim\Scal$ is $\C$-Mittag-Leffler if and only if it is a directed union of countably
presented submodules $N$ that belong to $\Mcal$ and satisfy $M/N\in\varinjlim\Scal$.
\end{Cor}
\begin{Proof}
By Corollary \ref{QML}, the module $M$  is $\C$-Mittag-Leffler if and only if it is a directed
union of countably presented submodules $N\in \varinjlim \Scal$ such that $N$ is
$\C$-Mittag-Leffler and the inclusion $N\subseteq M$ remains injective when tensoring with any
module $C \in \C$. The second condition means that $M/N\in\varinjlim\Scal$ by Lemma \ref{ft}(8).
The statement thus follows immediately from Theorem \ref{cotpair}.
\end{Proof}

An important example of cotorsion pairs of finite type is provided by tilting theory.

\medskip

\begin{Def} {\rm Let $n\in\N$.  A right $R$-module $T$ is {\em $n$-tilting} provided \\(T1) $T $ has projective dimension at most $n$,
\\(T2) Ext$^i_R(T, T^{(I)}) =0$ for each $i \geq 1$ and all sets $I$, and \\(T3) there
exist $r \geq 0$ and a long exact sequence \mbox{$0 \lora R \lora T_0 \lora \cdots \lora
T_r \lora 0$} such that $T_i \in $ Add$T$ for each $0 \leq i\leq r$.

\smallskip

\noindent
Cotorsion pairs  $(\Mcal, \Lcal)$ generated by some tilting module  $T$ are called {\em tilting cotorsion pairs}.
One then has $\Mcal\cap\Lcal=\Add T$.

\smallskip

\noindent
Dually, one defines   {\em cotilting} modules, and {\em cotilting cotorsion pairs}. If $(\C,\D)$ is a
 cotilting cotorsion pair cogenerated by a cotilting module $C$, then
$\C\cap\D=\Prod C$.
}
\end{Def}
\smallskip

Note that tilting cotorsion pairs are always of finite type, see \cite{bazher} and \cite{BaSt}.
Moreover, if $(\Mcal, \Lcal)$ is an $n$-tilting cotorsion pair generated by the tilting module $T$, then the cotorsion pair $(\C,\D)$
constructed as in Lemma \ref{ft} is an $n$-cotilting cotorsion pair  cogenerated by the cotilting module $T^\ast$, see \cite{AHT}.

\medskip

\begin{Cor}\label{tilting}
Let $(\Mcal,\Lcal)$ be a tilting cotorsion pair
  in $\Mod R$. Set
    $\Scal=\Mcal\cap\rmod R$, and   $\C=\Mcal^\intercal$.
    Let moreover $T$ be a tilting right $R$-module with $\Lcal=T^\perp$, and $C$ a  cotilting left $R$-module with $\C={}^\perp C$.
Then the following statements are equivalent for a right $R$-module  $M\in\varinjlim\Scal$.
\begin{itemize}
\item[(1)] $M$ is $\Lcal$-stationary.
\item[(2)] $M$ is $T$-stationary.
\item[(3)] $M$ is  $C$-Mittag-Leffler.
\item[(4)] $M$ is  $\C$-Mittag-Leffler.
\end{itemize}
Moreover, every $M\in\Mcal$ is strict $\Lcal$-stationary (and thus $\C$-Mittag-Leffler).  If $M$ is
countably generated, then $M\in\Mcal$ if and only if $M$ belongs to $\varinjlim\Scal$ and   is (strict)
$\Lcal$-stationary.
\end{Cor}
\begin{Proof}
Write $M$ as direct limit of a direct system $(F_\alpha,u_{\beta\alpha})$ of modules in $\Scal$.
Then for all $\alpha\in I$
 there is a short exact sequence $0\to F_\alpha\mapr{f_\alpha} B_\alpha\to B_\alpha/F_\alpha\to 0$   where $B_\alpha\in\Lcal$ and
 $B_\alpha/F_\alpha\in\Mcal$. Since $F_\alpha\in\Scal\subseteq\Mcal$, it follows that $B_\alpha\in\Lcal\cap\Mcal={\rm Add}T$. In particular,
 $f_\alpha$ is an $\Add T$-preenvelope.
  Further, $B_\alpha\in\Lcal\cap\varinjlim\Scal$, so  $B_\alpha^{\bullet}\in\C\cap\D={\rm Prod}C$ by \ref{ft}.
  Finally, since $B_\alpha/F_\alpha\in {}^\intercal\C$, the map $f\otimes C$ is a monomorphism for all $C\in\C$.  \\
  (1) and (4) are equivalent statements by Theorem \ref{cotpair}. Moreover, the implication (1)$\Rightarrow$(2) and (4)$\Rightarrow$(3) are trivial. \\
 (2)$\Rightarrow (4)$: Note that $M$ is Add$T$-stationary by Corollary \ref{closureB}.
 Now this implies (4) by \ref{Prop1}. \\ (3)$\Rightarrow (2)$: By \ref{definable}, $M$ is $Q$-Mittag-Leffler for all modules $Q$ that are pure submodules of products of copies of $C$. In particular,  we can take $Q=\bigoplus_{\alpha\in I} B_\alpha^{\bullet}$.
 The claim then follows from Proposition \ref{Prop2}.
\end{Proof}

\begin{Cor} \label{finitetilting}
  Let $T$ be a tilting module, $S=\End_R T$. Then the following hold true.\\
  (1)  Every finitely generated $S$-submodule of $T$ is a finite matrix
  subgroup.\\
  (2) ${}_ST$ is noetherian if and only if $T^\bullet$ is $\Sigma$-pure-injective.
  \end{Cor}
  \begin{Proof}
(1)  Let $(\mathcal M,\mathcal L)$ be the tilting cotorsion pair generated by $T$.  We know that  $T$ belongs to the kernel $\mathcal M\cap\mathcal L$. Thus  $T$ is strict $T$-stationary
  by Corollary~\ref{tilting}.
   By Theorem~\ref{dcc}, all matrix subgroups of $T$ of the form $H_{T,\,x}(T)$ with $x\in T$ are finite matrix subgroups, and
 of course, the  matrix subgroups of such form are precisely the cyclic $S$-submodules of $T$.
Observe that the class of finite matrix subgroups is closed under finite sums. So, we infer that
every finitely generated $S$-submodule of $T$, being
 a finite sum of cyclic submodules, is a finite matrix subgroup.\\
 (2) We know from  \cite[Proposition 3]{ZHZ} that $T^\bullet$ is $\Sigma$-pure-injective if and only if
$T$ satisfies the ascending chain condition on finite matrix subgroups. By (1), the latter means  that ${}_ST$ is noetherian.
\end{Proof}

If $(\mathcal M,\mathcal L)$ is a cotorsion pair of finite type,
then it follows from Theorem~\ref{cotpair} that $\Mcal$ is contained
in the class of strict $\Lcal$-stationary modules. If $\Scal =\Mcal
\cap \mathrm{mod}$-$R$, then the countably generated modules in
$\varinjlim \Scal$ that are strict $\Lcal$-stationary are precisely
the countably generated modules in $\Mcal$,  and they also coincide
with the countably generated modules in $\varinjlim \Scal$ that are
$\Lcal$-stationary. Raynaud and Gruson in \cite[p.~76]{RG} provide
examples showing that, in general, $\Mcal$ is properly contained in
the class of modules in $\varinjlim \Scal$ that are strict
$\Lcal$-stationary, and the latter class is  also properly contained
in the class of $\Lcal$-stationary modules. We explain  these
examples for  completeness' sake. First we prove the following

\begin{Lemma}\label{flatnotstrict}  Let $R$ be a ring, and let $F_1$ and $F_2$ be flat right $R$-modules such that there
exists an exact sequence
\[0\to R\stackrel{u}{\to} F_1\to F_2\to 0\]
\begin{itemize}
\item[(i)] Let $\Qcal$ be a class of left $R$-modules. If $F_2$ is $\Qcal$-Mittag-Leffler then so is $F_1$.
\item[(ii)] If $F_1$ is strict $R$-stationary, then $u$ splits.
\end{itemize}
\end{Lemma}

\begin{Proof} Statement $(i)$ follows from Examples~\ref{extp}(4).

$(ii)$ Assume that $F_1$ is strict $R$-stationary. Let $\Scal$ be
the class of finitely generated free modules. By
Proposition~\ref{strictdominatinglimits} and since $F_1\in
\varinjlim \Scal$, there exist $n>0$ and a homomorphism $v\colon
R\to R^n$, such that $u=tv$ for some $t\colon R^n\to F_1$ and
$H_u(R)=H_v(R)$. Since $u$ is a pure monomorphism  so is $v$, but
being a pure monomorphism between finitely generated projective
modules v  splits. Therefore the identity map belongs to
$H_v(R)=H_u(R)$, thus $u$ also splits.
\end{Proof}

The easiest instance of tilting cotorsion pair is the one generated by the class $\Scal$   of
finitely generated free modules. Then $\Scal ^\perp=\mathrm{Mod}$-$R$, and ${}^\perp(\Scal
^\perp)=\mathcal{P}$ is the class of all projective modules. The ring $R$ is a tilting module that
generates the tilting cotorsion pair $(\mathcal{P},\mathrm{Mod}$-$R )$. Note that   $\varinjlim
\Scal =\Fcal$ is the class of flat modules. The relative Mittag-Leffler,  strict  stationary
and stationary modules associated to this cotorsion pair in  Corollary~\ref{tilting} are the
Mittag-Leffler,  strict Mittag-Leffler and  $(\mathrm{Mod}$-$R)$-stationary modules, respectively.
If $F_1$ is a flat  Mittag-Leffler module and $\mathrm{Ext}_R^1(F_1,R)\neq 0$, then there is
a non split exact sequence
\[0\to R\to F_2\to F_1\to 0\]
In view of Lemma~\ref{flatnotstrict}, $F_2$ is a Mittag-Leffler module that is not strict
Mittag-Leffler. An   example for this situation is  the following:

\begin{Ex}\label{notstrict} For any set $I$, the abelian group $\mathbb{Z}^{I}$ is a flat strict Mittag-Leffler abelian group.
If $I$ is infinite, there exist nonsplit extensions of  $\mathbb{Z}^{I}$ by $\mathbb{Z}$, hence
there are  flat Mittag-Leffler abelian groups that are not strict Mittag-Leffler.
\end{Ex}

\begin{Proof} Of course, $\mathbb{Z}^{I}$ is flat. Since any finitely generated submodule of $\mathbb{Z}^{I}$ is contained in a
finitely generated direct summand of $\mathbb{Z}^{I}$ \cite[Proof of Theorem~19.2] {fuchs1}, we deduce from
 Proposition~\ref{strictdominatinglimits} that $\mathbb{Z}^{I}$ is   strict Mittag-Leffler.

When $I$ is infinite, $\mathbb{Z}^{I}$ is not a Whitehead group, that is,
$\mathrm{Ext}_\mathbb{Z}^1(\mathbb{Z}^{I}, \mathbb{Z})\neq 0$ \cite[Proposition~99.2]{fuchs2}.
Then the claim follows by the remarks above.
\end{Proof}


\section{The pure-semisimplicity conjecture}\label{pss}

  Throughout this section, we assume that  $R$ is a twosided  {artinian}, hereditary, indecomposable,
  left pure semisimple ring. It is well known that then every indecomposable finitely generated non-projective left $R$-module  is
   end-term of an almost split sequence in $R\LMod$
 consisting of finitely presented modules, and
 every    indecomposable finitely generated non-injective right $R$-module is the first term of an almost split sequence in $\Mod R$
  consisting of finitely presented modules.

  We adopt the notation $A=\tau\,C\,$ and
$C=\tau^{-}\,A$  if $\exs{A}{}{B}{}{C}$  is  an almost split sequence,  and define inductively $\tau^n$ resp.~$\tau^{-n}$.
We know from \cite{AV}
that there is a   {\em preprojective
component}  $\p$ in $\Mod R$, that is, a class of finitely generated indecomposable right $R$-modules
satisfying the following conditions.
\begin{enumerate}
\item For any $X\in\p$ there are a left almost split morphism $X\ra Z$ and a right almost split morphism $Y\to X$ in $\Mod R$ with $Z,Y$ being
finitely generated.
\item If $X\to Y$ is an irreducible map
with one of the modules lying in  $\p$, then both modules are in  $\p$.
\item  The Auslander-Reiten-quiver of $\p$ is connected and has no oriented cycles.
\item For every $Z\in \p$ there is $m\ge 0$ such that $\tau^m\,Z$ is projective.
\end{enumerate}
Similarly, there is a {\em preinjective component} in $R\LMod$, i.~e.~a class of finitely generated indecomposable left $R$-modules with the dual properties.
Moreover, the two components are related by the local duality, that is, there is a bijection  $ \q\rightarrow \p, \, _RA\mapsto A^+_{\,R}$.
The modules in add$\p$ are called \emph{preprojective}, the modules in add$\q$ are called \emph{preinjective}.

In \cite{key},
  the cotorsion pair $(\Mcal,\Lcal)$ in $\Mod R$  generated by $\p$ and the cotorsion pair $(\C,\D)$ in $R\LMod$  cogenerated by $\q$ are investigated.
  In particular, it is shown that there   is a finitely generated product-complete tilting and cotilting left $R$-module $W$  such that
  $\C =\text{Cogen}W={}^\perp W$ and $\D=\text{Gen}W=W^\perp$. Note that $\C=\Mcal^\intercal=\p^\intercal$ by Lemma \ref{ft}(2).
  Moreover, $(\Mcal,\Lcal)$ is a tilting cotorsion pair in $\Mod R$ with corresponding  cotilting cotorsion pair $(\C,\D)$ in $R\LMod$. But
  $(\C,\D)$ is also a tilting cotorsion pair in $R \LMod$, and the corresponding  cotilting cotorsion pair in $\Mod R$ is
 $(\varinjlim\add \p,\Ecal)$, see \cite[5.2 and 5.4]{key}.

 We now apply our previous results to this setup,     specializing to the case where $B^\bullet$ denotes the local dual of a module $B$. Let us  fix a tilting module $T$  such that $T^\perp=\Lcal$.

\begin{Prop}\label{know} The following statements hold true.
\begin{itemize}
\item[(1)]
$T$ is noetherian over its endomorphism ring.
\item[(2)]
 All right $R$-modules in $\Mcal$ are strict $T$-stationary (and hence $W$-Mittag-Leffler).
 \end{itemize}
 \end{Prop}
 \begin{Proof}
 (1)
Any left $R$-module is pure-injective, thus $\Sigma$-pure-injective.  In
particular,  $T^+$ is $\Sigma$-pure-injective. By Corollary \ref{finitetilting} we conclude that $T$ is noetherian over its endomorphism ring.
\\ (2) holds true by  Corollary \ref{tilting}.
 \end{Proof}

We remark that all left $R$-modules are Mittag-Leffler modules by \cite{AF},
hence strict Mittag-Leffler modules, see \ref{dualizing}(2) and \ref{firstremarks}.

As shown in \cite{key}, the validity of the Pure-Semisimplicity Conjecture is related to the question whether $W$ is endofinite. We obtain the following criteria for endofiniteness of $W$.
 \begin{Prop}\label{want}
 The following statements are equivalent.
\begin{itemize}
\item[(1)]
$W$ is endofinite.
\item[(2)]
$W^+$ is $W$-Mittag-Leffler.
\item[(3)]
 Every  (countable) direct limit of preprojective right $R$-modules is  $W$-Mittag-Leffler.
 \item[(4)]
 Every (countable)  direct system of preprojective right $R$-modules is $T$-stationary.
 \item[(5)]
 Every (countable) direct system of preprojective right $R$-modules has  limit in $\Mcal$.
 \item[(6)]
 If $A$ is direct limit of a (countable) direct system of preprojective right $R$-modules, and $L$ is a locally pure-injective module from $\Lcal$, then $\Ext^1_R(A,L)=0$.
 \end{itemize}
 \end{Prop}
 \begin{Proof}
 (1)$\Rightarrow$(5) follows from \cite[5.6]{key}, which asserts that $W$ is endofinite if and only if the class $\Mcal$ is closed under direct limits.\\
 (5)$\Leftrightarrow$(4)$\Leftrightarrow$(3) follows from Corollaries    \ref{reducetocountable} and
  \ref{tilting}.\\
 (3)$\Rightarrow$(2):  $W\in\C\cap\D$, hence $W^+\in\varinjlim\add \p$ by Lemma \ref{ft}(6).
  \\(2)$\Rightarrow$(1) holds by Proposition \ref{endofinite}.\\
 (5)$\Leftrightarrow$(6):  Of course, (5) implies (6). For the converse implication,
 observe first that $\Lcal$ is a definable class, hence it is closed under pure-injective envelopes.
  So, every module $L\in\Lcal$ is isomorphic to  a pure submodule of a  (locally) pure-injective module in  $\Lcal$.
   Note further that  the class
    $\Lcal'$  of all locally pure-injective modules from $\Lcal$ is closed under direct sums, because this is true for the tilting class $\Lcal$ and for the class of locally pure-injective modules, see \cite[2.4]{Zlpi}.
 Consider now a module $A$ which is direct limit of a countable direct system of preprojective right $R$-modules. Then   $A$ is countably presented,  and it follows from (6) and  \cite[2.5]{bazher} that  $A^\perp$ contains all pure submodules of modules in $\Lcal'$. We conclude that  $\Ext^1_R(A,L)=0$ for all $L\in \Lcal$, which proves $A\in\Mcal$.
  \end{Proof}

 \begin{Remark}
 {\rm
It is well known that every cotilting class is a torsion-free
class.  So, let us consider  the torsion pair defined by the
cotilting class $\varinjlim\add \p$, denoting by  $t$  the
corresponding torsion radical. Assume the following condition
holds true (cf. \cite[4.1]{L2}):
\begin{quote} If $N$ is a finitely generated submodule of $W^+$, then $t(W^+/N)$ is finitely generated.
  \end{quote}
Then it follows  that  $W$ is endofinite. In fact, since $W^+\in
\varinjlim\add \p$, all its finitely generated submodules  are in
$\Mcal$. Moreover, for every finitely generated submodule $N$ of
$W^+$ there is a finitely generated submodule $N'\subseteq W^+$
for which $N\subseteq N'$ and $W^+/N'\in\varinjlim\add \p$.
 To see this,
choose $N'$ such that $N'/N=t(W^+/N)$ and use that $N'/N$ is
finitely generated. So, we conclude that $W^+$ is a directed union
of finitely generated submodules  $N'$ that belong to $\Mcal$ and
satisfy $W^+/N'\in\varinjlim\add \p$.
 But then $W^+$ is $W$-Mittag-Leffler by   Corollary \ref{CML}, which means that $W$ is endofinite by  Proposition \ref{want}.
} \end{Remark}

We close with a  criterion for $R$ being of finite representation type.

\begin{Prop} Let $P$ be the direct sum of a
set of  representatives of the isomorphism classes of the modules
in $\p$. Then
 the following statements are equivalent.
\begin{enumerate}
\item $R$ has finite representation type.
\item Every countable direct system of preprojective modules is (strict)
$P$-stationary.
\end{enumerate}
\end{Prop}

\begin{Proof} It is clear that $(1)$ implies $(2)$. Indeed, $R$ is of
finite representation type if and only if all right (and left)
$R$-modules are Mittag-Leffler, hence stationary with respect to
any module \cite{AF}.

Assume $(2)$ holds.
To prove that $R$ has finite representation type,  it suffices to show that
$\p$ is finite, see  \cite[3.5]{key}. Assume on the contrary that $\p$ is infinite. By
applying \cite[Theorem 9]{ZHZ} to $\p$, we know that there are an
infinite family $(P_i)_{i\in \N}$ of  pairwise non-isomorphic
modules  in $\p$ and a sequence of homomorphisms $(f_i\colon
P_{i}\to P_{i+1})_{i\in \N}$ such that $f_i\dots f_0\neq 0$ for
any $i\in \N$.

Since all modules in $\p$ are endofinite by \cite[6.2]{AV}, the direct system
\[P_1\stackrel{f_1}\to P_2\stackrel{f_2}\to P_3\to\dots\to P_{n}\stackrel{f_n}\to P_{n+1}\to\dots\]
 is (strict) $\p$-stationary, see Corollary \ref{sigmapi}(3).
But by our assumption it is even  $P$-stationary, and since  $\Add P=\Add\p$, we infer from
 Corollary \ref{closureB}(3) that it is
$\Add\p$-stationary.
 By Corollary \ref{closureB}(2) it follows that for any $n\in \N$ there exists $m\ge n$
such that
\[ H_{f_{m}\cdots f_n}(P_{i})=
\bigcap _{l\ge n} H_{f_{l}\cdots f_n}(P_{i})\]
 for all $i\in \N$.
On the other hand, as our modules are preprojective, for any $i\in \N$ there exists $l_i$ such that
$\mathrm{Hom}_R(P_l, P_i)=0$ for any $l> l_i$, hence
$H_{f_l\cdots f_n}(P_i)=0$ for any $l\ge l_i$.
So, we deduce that $H_{f_{m}\cdots f_n}(P_{i})=0$ for all $i\in \N$.
In particular, we have
\[f_{m}\cdots f_n=\mathrm{Id}_{P_{m+1}}\, f_{m}\cdots f_n\in
H_{f_{m}\cdots f_n}(P_{m+1})=0,\] which contradicts the choice
of the sequence $(f_i)_{i\in \N}$. Thus we conclude that $\p$ is
finite.
\end{Proof}




\end{document}